\documentclass[12pt]{amsart}
\renewcommand{\bar}{\overline}
\usepackage{amsmath,amsthm,amscd,euscript}
\def \r{\mathbb R}

\def \z{\mathbb Z}

\newtheorem{theorem}{Theorem}[section]
\newtheorem{lemma}[theorem]{Lemma}

\newtheorem{proposition}[theorem]{Proposition}
\newtheorem{corollary}[theorem]{Corollary}
\theoremstyle{remark} \theoremstyle{definition}
\newtheorem{remark}[theorem]{Remark}
\newtheorem{definition}[theorem]{Definition}

\newtheorem{problem}{Problem}

\title[Completely empty pyramids on
integer lattices]{Completely empty pyramids on integer lattices
and two-dimensional faces of multidimensional continued
fractions.}
\author{O.~N.~Karpenkov}
\date{27 August 2006}
\thanks{Partially supported
by NWO-RFBR 047.011.2004.026 (RFBR 05-02-89000-NWO\_a) grant, by
RFBR SS-1972.2003.1 grant, by RFBR 05-01-02805-CNRSL\_a grant,
and by RFBR grant 05-01-01012a.}



\email[Oleg Karpenkov]{karpenk@mccme.ru}

\address{CEREMADE - UMR 7534 -- Universit\'e Paris-Dauphine,
France, 75775 Paris SEDEX 16}

\begin{document}
\input epsf

{\abstract In this paper we develop an integer-affine
classification of three-dimen\-sional multistory completely empty
convex marked pyramids. We apply it to obtain the complete lists
of compact two-dimensional faces of multidimensional continued
fractions lying in planes at integer distances to the origin equal
$2$, $3$, $4$, $\ldots$ The faces are considered up to the action
of the group of integer-linear transformations. In conclusion we
formulate some actual unsolved problems associated with the
generalizations for $n$-dimensional faces and more complicated
face configurations.}}

\maketitle \sloppy \normalsize

\tableofcontents

\section*{Introduction and background}

The main purpose of the present paper is to develop an
integer-affine classification of three-dimensional multistory
completely empty convex marked pyramids. We apply it to deduce an
integer-linear classification of compact two-dimensional faces of
multidimensional continued fractions in the sense of Klein lying
in planes at integer distances to the origin greater than unity.
Classification of two-dimensional faces leads to new algorithms
of two-dimensional continued fraction calculations. It is also
the first step in studying the combinatorial structure of
multidimensional continued fractions.

\subsection{General definitions}

 Consider a vector space $\r^{n+1}$ for some $n\ge
1$ over $\r$. A point or vector of $\r^{n+1}$ is called {\it
integer} if  all its coordinates are integers.

Consider some $k$-dimensional plane of $\r^{n+1}$. The
intersection of a finite number of closed $k$-dimensional
half-planes of the plane is said to be a {\it convex $($solid$)$
$k$-dimensional polyhedron} if it is homeomorphic to
$k$-dimensional closed disk. For $k=2$, $1$, or $0$ we have a {\it
convex polygon}, a {\it segment}, or a {\it point}. Here we
consider polyhedra as convex hulls (with all their inner points).

A convex polyhedron is said to be a {\it convex marked pyramid}
with some marked face and a vertex outside the plane containing
the face if it coincides with the union of all segments joining
the marked vertex with each point of the marked face. The marked
face is called the {\it base} of the marked convex pyramid and
the marked vertex --- the {\it vertex} of the marked convex
pyramid. A polyhedron is called a {\it convex pyramid} if some
structure of convex marked pyramid can be introduced for it.

A convex polyhedron (polygon, segment) is said to be {\it
integer} if all its vertices are all integer points. A convex
(marked) pyramid is said to be {\it integer} if it is an integer
convex polyhedron.

\begin{definition}
An integer convex polyhedron is called {\it empty} if it does not
contain integer points different from the vertices of the
polyhedron. An integer convex marked pyramid is called {\it
completely empty} if it does not contain integer points different
from the marked vertex and from the integer points of the base.
\end{definition}

Two sets in $\r ^{n+1}$ are said to be {\it integer-affine
equivalent} (or have the same {\it integer-affine type}), if
there exists an affine transformation of $\r^{n+1}$ preserving
the set of all integer points, and transforming the first set to
the second. Two sets in $\r ^{n+1}$ are said to be {\it
integer-linear equivalent} (or have the same {\it integer-linear
type}), if there exists a linear transformation of $\r^{n+1}$
preserving the set of all integer points, and transforming the
first set to the second.

\begin{definition}
A plane is called {\it integer} if it is integer-affine equivalent
to some plane passing through the origin and containing a
sublattice of the integer lattice, and the rank of the sublattice
is equivalent to the dimension of the plane.
\end{definition}

Consider some integer $(k-1)$-dimensional plane and an integer
point in the complement to this plane. Let the Euclidean distance
from the given point to the given plane equal $l$.  The minimal
value of nonzero Euclidean distances from all integer points of
the ($k$-dimensional) span of the the given plane and the given
point to the plane is denoted by $l_0$. Note that $l_0$ is always
greater than zero and can be obtained for some integer point of
the described span. The ratio $l/l_0$ is said to be the {\it
integer distance} from the given integer point to the given
integer plane.

\begin{figure}
$$
\epsfbox{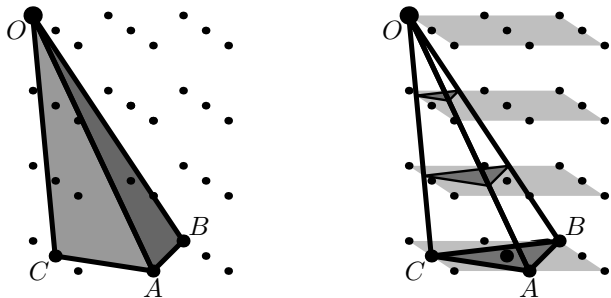}
$$
\caption{Two images of a completely empty three-story marked
pyramid with vertex $O$ and base $ABC$.}\label{pyramid}
\end{figure}

\begin{definition}
An integer convex marked pyramid is called {\it $l$-story} for
some positive integer $l$ if the integer distance from the vertex
of this pyramid to its base plane is equal to $l$. An integer
convex marked pyramid is called {\it multistory}/{\it
single-story} if the integer distance from the vertex of this
pyramid to its base plane is greater than one/equal to one. (See
example on Figure~\ref{pyramid}.)
\end{definition}

For any convex polygon there exist a single-story integer
three-dimensional convex marked pyramid with the given polygon as
the base (since any single-story integer convex marked pyramid is
completely empty). Two single-story three-dimensional convex
marked pyramids are integer-affine equivalent iff their bases are
integer-affine equivalent.

It turns out that the case of multistory convex marked pyramids
is essentially different from the single-story case. Only
polygons of a few integer-affine types can be bases of multistory
convex marked completely empty pyramids. For example, the
parallelogram with vertices $(0,0)$, $(0,1)$, $(1,1)$ and $(1,0)$
is not of that type. Besides, there exist integer-affine
nonequivalent multistory convex marked completely empty pyramids
whose bases are integer-affine equivalent.

In Section~1 of the present paper, we give the complete list of
integer-affine types of integer multistory convex marked
completely empty pyramids. To classify the pyramids, we study
arrangements of integer sublattices on the planes parallel to the
bases of the pyramids.

\subsection{Definition of multidimensional continued fra\-c\-ti\-ons in the
sense of Klein}

 The problem of generalizing
ordinary continued fractions to the higher-dimensional case was
posed by C.~Hermite~\cite{Herm} in 1839. A large number of
attempts to solve this problem lead to the birth of several
different remarkable theories of multidimensional continued
fractions. In this paper we consider the geometrical
generalization of ordinary continued fractions to the
multidimensional case presented by F.~Klein in~1895 and published
by him in~\cite{Kle1} and~\cite{Kle2}.

Consider a set of $n{+}1$ hyperplanes of $\r^{n+1}$ passing
through the origin in general position. The complement to the
union of these hyperplanes consists of $2^{n+1}$ open orthants.
Let us choose an arbitrary orthant.

\begin{definition}
The boundary of the convex hull of all integer points except the
origin in the closure of the orthant is called the {\it sail}.
The set of all $2^{n+1}$ sails of the space $\r^{n+1}$ is called
the {\it $n$-dimensional continued fraction} associated to the
given $n{+}1$ hyperplanes in general position in
$(n{+}1)$-dimensional space.
\end{definition}

Note that the union of all sails of any continued fraction is
centrally symmetric.

\begin{figure}
$$
\epsfbox{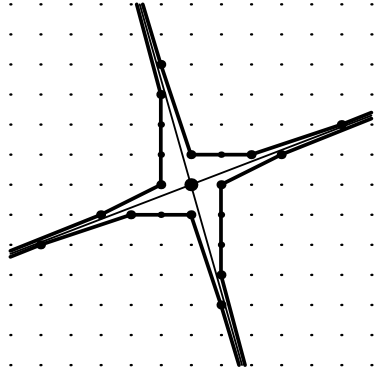}
$$
\caption{A one-dimensional continued fraction.}\label{cfrac}
\end{figure}

On Figure~\ref{cfrac} we show an example of one-di\-men\-si\-onal
continued fraction. This continued fraction contains four sails
(four broken lines on Picture~\ref{cfrac}). A description of
connections between ordinary continued fractions and geometrical
one-dimensional continued fractions can be found
in~\cite{KarItrig}, \cite{Hir}, and~\cite{Jun}.

Two $n$-dimensional continued fractions are said to be {\it equivalent}
if there exists a linear transformation that preserves the integer lattice
of the $(n+1)$-dimensional space and maps the sails of the first
continued fraction to the sails of the other.

Multidimensional continued fractions in the sense of Klein have
many relations with other branches of mathematics. For example,
J.-O.~Moussafir~\cite{Mou1} and O.~N.~German~\cite{Ger} studied
the connection between the sails of multidimensional continued
fractions and Hilbert bases. In~\cite{Tsu} H.~Tsuchihashi found
the relationship between periodic multidimensional continued
fractions and multidimensional cusp singularities, which
generalizes the relationship between ordinary continued fractions
and two-dimensional cusp singularities. M.~L.~Kontsevich and
Yu.~M.~Suhov discussed the statistical properties of the boundary
of a random multidimensional continued fraction in~\cite{Kon}.
The combinatorial topological generalization of Lagrange theorem
was obtained by E.~I.~Korkina in~\cite{Kor1} and its algebraic
generalization by G.~Lachaud~\cite{Lac}.

Theory of ordinary continued fractions was described, for example,
by A.~Ya.~Hinchin in~\cite{Hin}. V.~I.~Arnold presented a survey
of geometrical problems and theorems associated with
one-dimensional and multidimensional continued fractions in his
article~\cite{ArnPT} and his book~\cite{Arn2}). For the
algorithms of constructing multidimensional continued fractions,
see the papers of R.~Okazaki~\cite{Oka},
J.-O.~Moussafir~\cite{Mou2} and the author~\cite{Kar4}.

E.~Korkina in~\cite{Kor0},~\cite{Kor2},~\cite{Kor3}
and G.~Lachaud in~\cite{Lac},~\cite{Lac2},
A.~D.~Bruno and V.~I.~Parusnikov in~\cite{BP}, \cite{Par1},
\cite{Par1.1}, \cite{Par1.2} and~\cite{Par2}, the author in~\cite{Kar1}
and~\cite{Kar2} produced a large number of fundamental domains
for periodic algebraic two-dimensional continued fractions.
A nice collection of two-dimensional continued fractions
is given in the work~\cite{site} by K.~Briggs.

Besides the multidimensional continued fractions in the sense of Klein,
there exist several different generalizations of continued fractions
to the multidimensional case.
Some other well-known generalizations of continued fractions can
be found in the works of H.~Minkowski~\cite{Min}, G.~F.~Voronoi~\cite{Voro},
A.~K.~Mittal and A.~K.~Gupta~\cite{Mit1}, A.~D.~Bryuno and V.~I.~Parusnikov~\cite{BP2},
V.~Ya.~Skorobogat'ko~\cite{Sko}, and V.~I.~Shmoilov~\cite{Shm}.

\subsection{Two-dimensional faces of multidi\-men\-si\-o\-nal continued fractions}

Many classical papers were dedicated to studying algebraic and algorithmic properties
of multidimensional continued fractions.
The interest to geometrical properties of multidimensional continued fractions
was revived by V.~I.~Arnold's work~\cite{Arn1}
and the subsequent paper of E.~I.~Korkina~\cite{Kor0} on the classification of
$A$-algebras with three generators.
In~1989 and later, V.~I.~Arnold formulated a series of problems and conjectures
associated to the geometrical and topological properties of sails
of multidimensional continued fractions. The majority of these problems
are still open. The geometry of sails has not been sufficiently studied.

In the present work, we make the first steps in the investigation of
geometric properties of sails.
One of the first natural questions here is the following:
{\it what compact faces can sails of multidimensional continued fractions have?
$($these objects are usually studied up to the integer-linear equivalence relation$)$?}

The complete answer to this question was known only for
one-dimensional continued fractions.
{\it For any non-negative integer number $n$ there exists a one-dimensional face
with exactly $n$ integer points inside.
Two compact faces with the same numbers of integer points inside are integer-linear equivalent.}

In the two-dimensional case the original question decomposes into two questions.

{\it  What compact faces contained in planes at integer distances
from the origin equal to one can sails of multidimensional
continued fractions have $($up to integer-linear equivalence$)$?}

{\it What compact faces contained in planes at integer distances
from the origin greater than one can sails of multidimensional
continued fractions have $($up to integer-linear equivalence$)$?}

The answer to the first question is quite straightforward. For
any convex polygon $P$ at the unit integer distance from the
origin, there exist an integer positive $k$, a $k$-dimensional
continued fraction, and some face $F$ of a sail of this continued
fraction, such that $F$ is integer-affine equivalent to $P$.
Furthermore, two two-dimensional faces in the planes at the unit
integer distance from the origin are integer-linear equivalent iff
the corresponding polygons are integer-affine equivalent.

Note that up to this moment the following statement on compact
two-dimensional faces (of sails of multidimensional continued
fractions) contained in planes at integer distances from the
origin greater than one was known. {\it Such faces are either
triangles or quadrangles} (see the work~[3] by J.-O.~Moussafir).

In the present work we classify compact two-dimensional faces
contained in planes at integer distances from the origin greater
than one up to integer-linear equivalence. This result was
announced in~\cite{Kar5}. We give the complete lists for continued
fractions of any dimension. This result is based on the
classification of three-dimensional multistory completely empty
convex marked pyramids.

\subsection{Description of the paper}

We start in Section~1 with introducing Theorem~A  on
integer-affine classification of three-dimensional multistory
completely empty convex marked pyramids. In this section we also
formulate Theorem~B on integer-linear classification of
two-dimensional faces of the sails at integer distance greater
than one. The integer-affine classification of two-dimensional
faces contained in planes at integer distances from the origin
greater than one (Corollary~C) directly follows from the
integer-linear classification of two-dimensional faces contained
in planes at integer distances from the origin greater than one.
In Sections~2 and~3 we prove Theorem~A and Theorem~B respectively.
And, finally, in Section~4 we give a list of unsolved problems
associated with main theorems of this work.

\section{Formulation of main results}

\subsection{Classification of two-dimensional multistory completely empty pyramids}

By $(a_1, \ldots , a_k)$ in $\r ^n$ for $k<n$ we denote the point
$(a_1, \ldots,a_k, 0, \ldots , 0)$.

Denote the marked pyramid with vertex at the origin and
quadrangular base $(2,-1,0)$, $(2,-a-1,1)$, $(2,-1,2)$, $(2,b-1,1)$,
where $b\ge a \ge 1$, by $M_{a,b}$.

Denote the marked pyramid with vertex at the origin and
triangular base \\
$(\xi,r-1,-r)$, $(a+\xi,r-1,-r)$, $(\xi,r,-r)$,
where $a \ge 1$, $r\ge 1$, by $T_{a,r}^{\xi}$;\\
$(2,1,b-1)$, $(2,2,-1)$, $(2,0,-1)$, where $b\ge 1$, by $U_b$;\\
$(2,-2,1)$, $(2,-1,-1)$, $(2,1,2)$ by $V$;\\
$(3,0,2)$, $(3,1,1)$, $(3,2,3)$ by $W$ (pyramid $W$ is shown on
Figure~\ref{pyramid}).

The integer-affine types of the bases of the described above triangular
and quadrangular pyramids are shown on Figure~1.

\begin{figure}
$$
\epsfbox{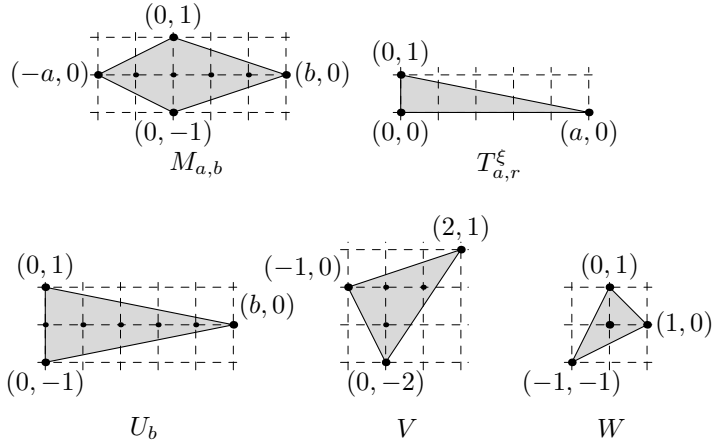}
$$
\caption{The integer-affine types of the bases of the marked pyramids
of List ``M-W''.}\label{osnovaniya}
\end{figure}

\vspace{1mm} {\bf Theorem A.} {\it Any multistory completely
empty convex three-dimensional marked py\-ra\-mid is
integer-affine equivalent exactly to one of the marked pyra\-mids
from the following list.

{\bf List ``M-W'':}\label{sp_M-W}

\nopagebreak

--- the quadrangular marked pyramids $M_{a,b}$, with integers $b\ge a \ge 1$;

--- the triangular marked pyramids $T_{a,r,}^{\xi}$,
where $a \ge 1$, and $\xi$ and $r$ are relatively prime, and
$r\ge 2$ and $0<\xi\le r/2$;

--- the triangular marked pyramids $U_b$, where $b\ge 1$;

--- the triangular marked pyramid $V$;

--- the triangular marked pyramid $W$.
} \vspace{1mm}

We give the proof of Theorem~A in Section~\ref{pA}.

\subsection{Compact two-dimensional faces at distance greater than one}

Note that up to this moment the following statement on compact
two-dimensional faces contained in planes at the integer distance
from the origin greater than one was known.

\vspace{1mm} {\bf Theorem (J.-O.~Moussafir~\cite{Mou2}.)} {\it Let
$F$ be a two-dimensional compact face of some sail of a
two-dimensional continued fraction. Let $r$ be the integer
distance from the origin to the plane, containing the face.

1. If $r=1$, $F$ may have arbitrary many vertices.

2. If $r=2$, $F$ has at most 4 vertices.

3. If $r\ge 3$, $F$ has three vertices. \qed
}\vspace{1mm}

Here we present a new theorem on integer-linear classification
and its corollary on integer-affine classification of
two-dimensional faces of multidimensional sails (the faces are
again contained in the planes at integer distances greater than
one from the origin). Note that from this theorem and its
corollary it follows that the second item of Moussafir's theorem
can be strengthened:

{\it 2 $'$. If $r=2$, $F$ has three vertices.}

Quadrangular faces for the case of $r=2$ are possible only for
$n$-dimensional continued fractions where $n\ge 3$.

\vspace{1mm} {\it {\bf Theorem~B.}
Any compact two-dimensional face of a sails of a two-dimensional
continued fraction contained in a plane at an integer distance
from the origin greater than one is integer-linear equivalent
exactly to one of the faces of the following list.

{\bf List ``$\alpha_{2}$'':}\label{sp_alpha}

--- triangle with vertices $(\xi,r-1,-r)$, $(a+\xi,r-1,-r)$,
$(\xi,r,-r)$, where $a \ge 1$, integers $\xi$ and $r$ are
relatively prime and satisfy the following inequalities $r\ge 2$
and $0<\xi\le r/2$;

--- triangle with vertices $(2,1,b{-}1)$, $(2,2,-1)$, and $(2,0,-1)$ for $b\ge 1$;

--- triangle with vertices $(2,-2,1)$, $(2,-1,-1)$, and  $(2,1,2)$;

--- triangle with vertices $(3,0,2)$, $(3,1,1)$, and $(3,2,3)$.\\
All triangular faces of List ``$\alpha_2$'' are realizable by
sails of dimension two and integer-linear nonequivalent to each
other.

Any compact two-dimensional face of a sails of a $n$-dimensional
$(n\ge 3)$ continued fraction contained in a plane at an integer
distance from the origin greater than one is integer-linear
equivalent exactly to one of the faces of the following list.

{\bf List ``$\alpha_n$'', $n\ge 3$:}

--- quadrangle with vertices $(2,-1,0)$, $(2,-a-1,1)$, $(2,-1,2)$,
$(2,b-1,1)$, where $b\ge a \ge 1$,

--- triangle with vertices $(\xi,r-1,-r)$, $(a+\xi,r-1,-r)$,
$(\xi,r,-r)$, where $a \ge 1$, integers $\xi$ and $r$ are
relatively prime and satisfy the following inequalities $r\ge 2$
and $0<\xi\le r/2$;

--- triangle with vertices $(2,1,b-1)$, $(2,2,-1)$, and $(2,0,-1)$ for $b\ge 1$;

--- triangle with vertices  $(2,-2,1)$, $(2,1,2)$, and  $(2,-1,-1)$;

--- triangle with vertices $(3,0,2)$, $(3,1,1)$, and $(3,2,3)$.\\
All faces of List ``$\alpha_n$'' are realizable by sails of any
dimension greater than two and integer-linear nonequivalent to
each other. } \vspace{1mm}

\begin{remark}
Note that for any compact face of a sail we can associate an
integer completely empty convex marked pyramid with marked vertex
at the origin and this face as base. Therefore integer-affine
types of such marked pyramids are in one-to-one correspondence
with integer-linear types of faces (see lemma~\ref{l81} below).
\end{remark}

We give the proof of Theorem~B in Section~\ref{pB}.

{\bf Corollary~C.} {\it
Any compact two-dimensional face of a sails of a multidimensional
continued fraction contained in a plane at integer distance from
the origin equals $r$ is integer-affine equivalent exactly to one
of the polygons of the list $\beta_r$ shown below.

\begin{figure}
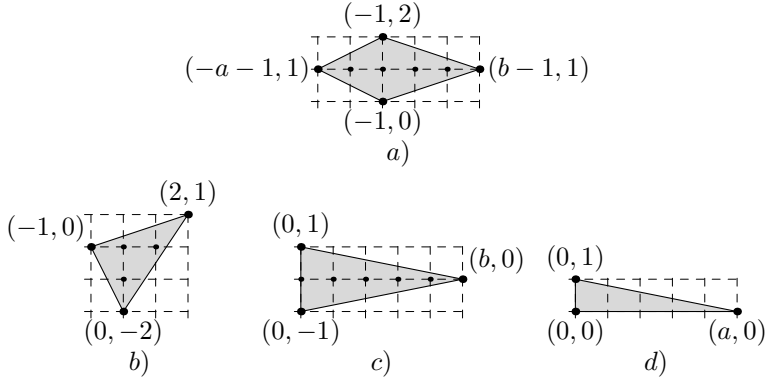

$$
\begin{array}{c}
\epsfbox{theorem1.11}\\
\epsfbox{theorem1.12}\\
\end{array}
$$
\caption{Integer-affine types of faces of List
``$\beta_2$''.}\label{beta2}
\end{figure}%

{\bf List ``$\beta_2$'':}\label{sp_beta}

\nopagebreak

--- quadrangle with vertices $(-1,0)$, $(-a-1,1)$, $(-1,2)$,
$(b-1,1)$, where $b\ge a \ge 1$ $($see the case of $a=2$, $b=3$
on Figure~\ref{beta2}a$)$$)$; quadrangular faces are possible
only for $n$-dimensional continued fractions where $n\ge 3$;

--- single triangle $(-1,0)$, $(0,-2)$,
$(2,1)$ $($see Figure~\ref{beta2}b$)$$)$;

--- triangle $(0,-1)$, $(0,1)$, $(b,0)$, for $b\ge 1$
$($see the case of $b=5$ on Figure~\ref{beta2}c$)$$)$;

--- triangle
$(0,0)$, $(a,0)$, $(0,1)$,  for  $a\ge 1$ $($see the case of $a=5$
on Figure~\ref{beta2}d$)$$)$.

\begin{figure}[ht]
$$\epsfbox{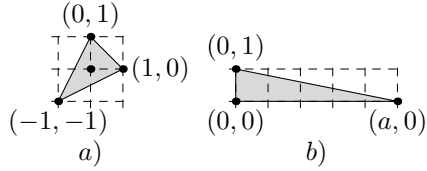}$$
\caption{Integer-affine types of faces of List
``$\beta_3$''.}\label{beta3}
\end{figure}%

{\bf List ``$\beta_3$'':}

\nopagebreak

--- single triangle $(-1,-1)$, $(1,0)$,
$(0,1)$ $($see Figure~\ref{beta3}a$)$$)$;

--- triangle
$(0,0)$, $(a,0)$, $(0,1)$,  for $a\ge 1$ $($see the case of $a=5$
on Figure~\ref{beta3}b$)$$)$.

\begin{figure}
$$\epsfbox{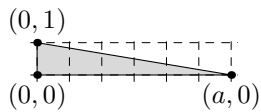}$$
\caption{Integer-affine types of faces of List ``$\beta_r$'', for
$r \ge 4$.}\label{betar}
\end{figure}

{\bf List ``$\beta_r$'', ($r\ge 3$):}

\nopagebreak

--- triangle with vertices $(0,0)$, $(a,0)$, and $(0,1)$,  for some $a\ge 1$
$($see the case of $a=6$ on Figure~\ref{betar}$)$, the
corresponding convex marked pyramid is integer-affine equivalent
to~$T_{a,r}^\xi$, where the integers $\xi$ and $r$ are relatively
prime and satisfy $0<\xi\le r/2$. For different $\xi$ the
corresponding faces are integer-linear nonequivalent but
integer-affine equivalent.

For any integer $r$ the faces of List $\beta_r$ are
integer-affine nonequivalent to each other; List $\beta_r$ is
irredundant. \qed }

\vspace{1mm}

The integer-affine and the integer-linear classifications
coincide, for $r<5$. For $r\ge 5$, the integer-linear
classification contains the integer-affine classification.

For any integers $n\ge 3$ and $r\ge 2$, the integer-linear
classification of compact two-dimensional faces contained in
planes at integer distances from the origin greater than one of
sails of $n$-dimensional continued fractions coincides with the
integer-affine classification of completely empty $r$-story
three-dimensional convex marked pyramids.

\section{Proof of Theorem~A}\label{pA}

\subsection{Preliminary definitions and statements}

Before proving the main theorem, we give several definitions and fix the notation,
and also formulate some general statements
that we will further use in the proof of the main statements.

For an integer polygon in some two-dimensional subspace
the ratio of its Euclidean volume
to the minimal possible Euclidean volume of an integer triangle in the same
two-dimensional subspace is called the{\it integer volume} of this polygon.

\begin{remark}
Our integer volume is a positive integer
(for a parallelogram, the usual volume will be two times less).
The integer volume of a triangle is equal to the index of the lattice
generated by its edges.
\end{remark}

An integer polyhedron (polygon) is called {\it empty},
if it does not contain integer points in its interior, and the set of integer points
of the faces coincides with the set of vertices of the polyhedron (polygon).

Let $ABCD$ be a tetrahedron with an ordered set of vertices $A$, $B$, $C$
and $D$. Denote by $P(ABCD)$ the following parallelepiped:
$$
\{A+ \alpha \bar{AB}+ \beta \bar{AC} +\gamma \bar{AD}|
\makebox[.2cm]{} 0\le\alpha\le1, 0\le\beta\le1, 0\le\gamma\le1\}.
$$

\begin{definition}
Now we specify some useful coordinates in the three-dimensional subspace
containing $P(ABCD)$ of $\r ^n$.
Let $b$, $c$, and $d$ be the distances from $B$, $C$, and $D$
to the two-dimensional planes containing the faces $ACD$, $ABD$, and $ACD$
respectively.
Let us define the coordinates of $A$, $B$, $C$, and $D$
as follows: $(0,0,0)$, $(b,0,0)$, $(0,c,0)$, and $(0,0,d)$ respectively.
The coordinates of all other points in this three-dimensional subspace are
uniquely defined by means of linearity. We call them the {\it integer-distance
coordinates} with respect to $P(ABCD)$.
\end{definition}

\begin{remark}
For any set of vertices $A$, $B$, $C$, and $D$ with the order as in $P(ABCD)$,
the integer-distance coordinates are uniquely defined.
\end{remark}

Using integer-distance coordinates by {\it integer lattice nodes
of $\r ^n$} (or, for short, {\it lattice nodes})
we mean integer points in the original coordinates in $\r ^n$.

\begin{remark}
Note that any lattice node of the three-dimensional space described above
has integer coordinates in the new integer-distance system of coordinates.
The inverse is not true. There exists an integer-distance system of coordinates
and such a point in the corresponding three-dimensional space
with integer coordinates which is not a lattice node.
For lattice nodes, the absolute values of their new coordinates coincide with the
integer distances from these lattice nodes to the planes containing the corresponding faces
of the parallelepiped.
\end{remark}

Let us continue with the following definition.

\begin{definition}
Two points $P$ and $Q$ are said to be
{\it equivalent with respect to some integer parallelogram} $ABCD$,
if there exist such integers $\lambda$ and $\beta$
that $P=Q+\lambda \bar{AB}+\beta \bar{AC}$.
The set of all equivalence classes of the integer lattice
with respect to the integer parallelogram $ABCD$
is called the {\it quotient-lattice} of the space by this integer parallelogram.
\end{definition}

Note that any equivalence class is contained in one of two-dimen\-sio\-nal planes
parallel to the plane of the parallelogram.

\begin{proposition}\label{m_st}
Consider an integer parallelepiped  $ABCDA'B'C'D'$ in $\r^3$
and some integer plane $\pi$ parallel to the face $ABCD$.
Let $\pi$ intersect the parallelepiped $($along a parallelogram$)$.
Then the following two statements hold.

First, $\pi$ contains only finitely many equivalence classes of the integer lattice
with respect to the integer parallelogram $ABCD$.
Their number is equivalent to the index of the sublattice generated by the edges
$\bar {AB}$ and $\bar{AC}$ in the integer lattice of the plane containing
$GBCD$.

Secondly, for any yquivalence class of the integer lattice contained in $\pi$
with respect to the integer parallelogram $ABCD$ exactly one of the following conditions holds.
\\
a$)$ only one point of the equivalence class is in the parallelogram,
it is an inner point of the parallelogram;
\\
b$)$ two points of the equivalence class are in the parallelogram,
they are contained in $($open$)$ opposite edges of the parallelogram;
\\
c$)$ four points of the equivalence class are in the
parallelogram, they coincide with vertices of the parallelogram.
\end{proposition}

We skip the proof of Proposition~\ref{m_st}.
It is straightforward and is based on the following easy lemma.

\begin{lemma}\label{l1}
Consider an integer parallelepiped with an empty face.
Let some parallel to this face plane intersect the parallelepiped $($along a parallelogram$)$.
Then exactly one the following statements holds.
\\
a$)$ only one integer point is in the parallelogram,
it is an inner point of the parallelogram;
\\
b$)$ two integer points are in the parallelogram,
they are contained in $($open$)$ opposite edges of the parallelogram;
\\
c$)$ four integer points are in the parallelogram, they coincide
with vertices of the parallelogram. \qed
\end{lemma}

\subsection{First results on empty integer tetrahedra}

In this subsection we present the corollary of White's
theorem~\cite{Whi} (see also~\cite{Ger}). Here without lose of
generality we consider only the three-dimensional space.

\begin{theorem}{\bf (G.~K.~White, 1964~\cite{Whi}.)}
Let $\Delta \subset \r ^3$ be an integer three-dimensional
simplex, let $E_i=\{ \sigma_i, \sigma_i '\}$, $i=1,2,3$ be a set
of points of two opposite edges $ \sigma_i, \sigma_i '$ for
$\Delta$. Then $(\Delta \setminus E_i)\cap \z ^3$ is empty iff
there exist a $j$ and two neighboring planes $\pi_j$, $\pi_j '$
$($by neighbor we mean that there is no integer points
``between'' these planes $\pi_j$ and $\pi_j '$$)$, such that
$\sigma_i \subset \pi_j$ and $\sigma_i ' \subset \pi_j '$. \qed
\end{theorem}

We will use the following corollary on empty integer tetrahedra
for the classification of empty convex multistory tetrahedra and
also further in the proof of Theorem~A.

\begin{corollary}\label{symplex}
Let $ADBA'$ be some empty integer tetrahedron.
Then all integer points of the parallelepiped $P(ADBA')$
are in the plane passing trough
two centrally-symmetric edges of the parallelepiped.
This two edges are not the edges of the tetrahedron $ADBA'$.
\qed
\end{corollary}

\begin{remark}
The number of planes passing through two centrally-symmetric edges of the parallelepiped
equals six, but only three of them do not contain the edges of the tetrahedron.
\end{remark}

For the proofs see~\cite{Whi}.

\subsubsection{Classification of empty triangular marked
pyra\-mids}

Corollary~\ref{symplex} allows to describe all integer-affine
types of empty triangular marked pyramids (i.e. tetrahedra with
one marked vertex each).

Let $r$ be some positive integer, and $\xi$ be nonnegative integer.
Denote by $P_{r}^{\xi}$ the marked pyramid with vertex at $(0,0,0)$
and the triangular base $(0,1,0)$, $(1,0,0)$, $(\xi,r-\xi,r)$.

\begin{corollary}\label{col_pyr}
Any integer  empty triangular marked pyramid is inte\-ger-affine equivalent to exactly
one of the pyramids of the

{\bf List ``P'':}

--- $P_1^0$;

--- $P_r^\xi$, where $\xi$ and $r$ are relatively prime, and $r{\ge} 2$, and $0{<}\xi{\le} r/2$.

All triangular marked pyramids of List~``P'' are empty and
integer-affine nonequivalent to each other.
\end{corollary}

\begin{proof}

{\bf 1. Completeness of List ``P''.}
Let us show that an arbitrary empty integer marked pyramid $ADBA'$
(with a vertex $A$) is integer-affine equivalent to one of the marked pyramids of ``P''.

Suppose that, the integer distance from its marked vertex to the plane containing the marked base
equals some positive integer $r$.
If $r=1$, then the vertices of the marked pyramid generate the three-dimensional
integer lattice, and therefore such a marked pyramid is integer-affine equivalent to $P_1^0$
(here $A$ corresponds to the marked vertex of $P_1^0$).

Suppose now that $r>1$. By Corollary~\ref{symplex},
all integer points of the parallelepiped $P(ADBA')$ are contained
exactly in one of the three planes passing through centrally-symmetric
edges of the parallelepiped and not containing the edges of the tetrahedron $ADBA'$.
Denote the vertices of the marked base $DBA'$ by $\bar B$,
$\bar D$, and $\bar A'$ in a such way that all inner integer points of the parallelepiped
$P(A \bar D \bar B \bar A')$ are contained in the plane passing
through $\bar B\bar D$ and the centrally-symmetric edge.

Consider the integer-distance coordinates with respect to the parallelepiped
$P(A \bar D \bar B \bar A')$.
Take the intersection of the parallelepiped with the plane $z=1$ in these coordinates.
There is only one lattice node in the intersection, by Corollary~\ref{symplex}
its coordinates are $(r-\xi,\xi,1)$.
Denote this lattice node by $K$.

If the integers $\xi$ and $r$ have some common integer divisor $c\ge 1$,
then the point with the coordinates $(\frac{r-\xi}{c}r,\frac{\xi}{c}r,c)$
is a lattice node.
Hence the point $(0,0,c)$ is also a lattice node.
And then the marked pyramid $A \bar D \bar B \bar A'$ is not empty.
Thus the integers $\xi$ and $r$ are relatively prime.

Since the integer distance from $K$ to the two-dimensional plane
containing the face $A \bar D \bar B$
equals one,
there exists an integer-affine transformation taking the tetrahedron $A\bar B\bar D K$
to the tetrahedron with vertices $(0,0,0)$,
$(0,1,0)$, $(1,0,0)$, and $(1,1,1)$.
Here the point $\bar A'$ maps to $(\xi,r-\xi,r)$.
Hence the integer-affine type of the marked pyramid $ABDA'$
coincides with the integer-affine type of the marked pyramid $A \bar B \bar D \bar A'$,
and it in turn coincides with the integer-affine type of the marked pyramid
$P_r^\xi$, where $0<\xi< r$, and $\xi$ and $r$ are relatively prime.
It remains to say that the marked pyramids $P_r^\xi$ and $P_r^{r-\xi}$
can be mapped one to another by the symmetry about the plane $x=y$
(which preserves the integer lattice).
Therefore the marked pyramids $P_r^\xi$ and $P_r^{r-\xi}$
are integer-affine equivalent.

\medskip

{\bf 2. Emptiness of the marked pyramids of List ``P''.}
Let us show that all listed marked pyramids $P_\xi^r$ are empty.

The intersection of the plane $z=b$ (for $1\le b \le (r-1)$) and marked pyramid $P_\xi^r$
is the triangle $A_kB_kD_k$ with the following coordinates of the vertices:
$$
\left( \frac{b}{r}\xi,\frac{b}{r}(r{-}\xi),b \right), \quad
\left( \frac{b}{r},\frac{b}{r}(r{-}\xi){+}\frac{r{-}b}{r},b
\right), \quad \left(
\frac{b}{r}\xi{+}\frac{r{-}b}{r},\frac{b}{r}(r{-}a),b \right).
$$
The triangle $A_kB_kD_k$ is contained in the band $b\le x+y \le b+\frac{r-b}{r}$, $z=b$.
This band contains only integer points with coordinates $(t,b-t,b)$ for integer $t$.
Hence it remains to check if $A_k$ is integer.
Since $\xi$ and $r$ are relatively prime and $d<r$, the first coordinate of
$A_k$ is not integer.
Therefore all marked pyramids $P_\xi^r$ of List ``P'' are empty.

\medskip

{\bf 3. Irredundance of List ``P''.}
We will show now that all marked pyramids $P_\xi^r$ of List ``P'' are integer-affine
nonequivalent to each other.
Note that the integer distance from the marked vertex to the plane containing the base
is an integer-affine invariant. Therefore the pyramids with nonequivalent parameter $r$
are integer-affine nonequivalent.

To distinguish the marked pyramids with the same $r$,
we construct the following integer-affine invariant.
Consider an arbitrary empty mar\-ked pyramid $ABDA'$  with a marked vertex $A$
and the corresponding trihedral angle also with vertex $A$ and triangle
$DBA'$ as its base.
By White's theorem, exactly one lattice node of the trihedral angle (we denote this lattice node by $K$)
is contained in the two-dimensional plane parallel to the face $DBA'$
and at the integer distance $r+1$ from~$A$.
By Corollary~\ref{symplex}, the integer distances from $K$
to two-dimensional planes of the angle are equal to $1$, $\xi$, $r-\xi$
(for some integer $\xi$).
The trihedral angle and $K$ are uniquely defined by the marked pyramid
up to the symmetries of the marked pyramid preserving the marked vertex.
The group of such symmetries permutes all integer distances
from $K$ to the planes containing the faces of the angle.
Hence, the unordered system of integers $[1, \xi, r-\xi]$
is an invariant.
This invariant distinguishes all marked pyramids $P_\xi^r$
with the same integer distance $r$.
\end{proof}

\begin{proposition}
Let relatively prime integers $\xi$ and $r$ satisfy the following
inequalities: $r\ge 2$, $0<\xi\le r/2$.
Then the marked pyramid $P_r^{\xi}$ is integer-affine equivalent
to the marked pyramid $T_{1,r}^{\xi}$.
\end{proposition}
\begin{proof}
The marked pyramid  $T_{1,r}^{\xi}$ is the image of $P_r^{\xi}$ under
the integer-linear transformation
$$
\left(
\begin{array}{ccc}
\xi +1 & \xi & -\xi \\
r-1 & r-1 & 2-r \\
-r & -r & r-1 \\
\end{array}
\right)
.
$$
\end{proof}

\begin{corollary}\label{col_pyr2}
Any integer empty $r$-story $($$r\ge 2$$)$ triangular marked pyramid
is integer-affine equivalent exactly to one of the marked pyramids
$T_{1,r}^\xi$ for relatively prime integers $\xi$ and $r$ satisfying $0<\xi\le r/2$.
All such pyramids $T_{1,r}^\xi$ are empty
$($and integer-affine nonequivalent if the corresponding parameters $r$ and $\xi$
do not coincide$)$.
\qed
\end{corollary}

\subsubsection{Classification of integer empty tetrahedra}

A certain difference between the integer-affine classification of
integer empty triangular marked pyramids (with marked vertex) and
the integer-affine classification of integer empty tetrahedra
(without mar\-ked vertices) occurs. The first steps in the
integer-affine classifications of integer empty tetrahedra were
made by J.-O.~Moussafir in~\cite{Mou2}.

\begin{theorem}{\bf (J.-O.~Moussafir~\cite{Mou2}.)}
Any integer empty tetrahedron is integer-affine equivalent
to the tetrahedron with vertices $(0,0,0)$,
$(1,0,0)$, $(0,1,0)$, and $(u,v,d)$, for some integers $u$, $v$ and $d$,
where $u$, $v$ and $u+v-1$ are relatively prime with $d$,
and one of the integers $u+v$, $u-1$, $v-1$ is divisible by $d$.
$($These tetrahedra are sometimes called
{\it Hermitian normal forms of the simplexes}.$)$
\end{theorem}

Note that many of such Hermitian normal forms are integer-affine equivalent
to each other.
The following consequence of Corollary~\ref{symplex} improves
Moussafir's theorem.

\begin{corollary}\label{col_symp}
Any integer empty tetrahedron is in\-te\-ger-affine equivalent
exactly to one of the following tetrahedra:

--- $P_1^0$;

--- $P_r^\xi$, where $r \ge 2$, $0<\xi< r$,
and the element $(\xi \mod r)$ of the additive group $\z/m\z$
is also contained in the associated multiplicative group $(\z/m\z)^*$
$($i.e. integers $\xi$ and $r$ are relatively prime$)$.

All listed integer tetrahedra are empty.
Two tetrahedra $P_{r_1}^\xi$ and $P_{r_2}^{\nu}$ are integer-affine
equivalent iff $r_1=r_2$ and $($for $r_1>1$$)$ one of the following
equalities in $(\z/m\z)^*$ holds:
$$
(\xi \mod r_1)=(\pm 1)\cdot (\nu \mod r_1)^{\pm 1}.
$$
\end{corollary}

\begin{proof}

{\bf 1. Completeness of the list.}
By Corollary~\ref{col_pyr}, any empty integer tetrahedron is
integer-affine equivalent to some tetrahedron of the list of
Corollary~\ref{col_symp}.

\medskip

{\bf 2. Emptiness of the tetrahedra of the list.}
By Corollary~\ref{col_pyr}, the tetrahedron $P_\xi^{r}$
is empty for relatively prime integers $r$ and $\xi$ satisfying $r \ge 2$ and $\xi\le r/2$.
Since $P_\xi^{r}$ and $P_\xi^{r-\xi}$ are integer-affine equivalent
and $P_1^0$ is empty, all tetrahedra of the list of Corollary~\ref{col_symp} are empty.

\medskip

{\bf 3. Proof of the last statement of Corollary~\ref{col_symp}.}
The (affine) symmetry group of the right tetrahedron ($S_4$) includes
the (affine) symmetry group of the right tetrahedron with marked vertex ($S_3$).
Now the list of the four trihedral angles associated with all four vertices
of the tetrahedron is uniquely defined.
We chose one lattice node for each of these angles as we did in the proof of the previous corollary.
Direct calculations show that the integer distances from these points to the
four two-dimensional planes containing the faces of the tetrahedron are
$$
(1,1,\xi, r-\xi), \quad (1,1,\xi,r-\xi), \quad (\nu,r-\nu,1,1),
\quad \mbox{and} \quad (\nu,r-\nu,1,1),
$$
where $(\xi \mod r)\cdot (\nu \mod r)=1$ in $(\z/m\z)^*$. The
set of these numbers up to the group $S^4$ of permutations action
(for all points at the same time) is an integer-affine invariant.
Therefore, the tetrahedra $P_{r_1}^\xi$, $P_{r_2}^\nu$,
$P_r^{r-\xi}$, and $P_r^{r-\nu}$ are integer-affine equivalent
and the invariant distinguishes all other tetrahedra.
\end{proof}

\begin{remark} The integer-affine classifications of integer empty
triangular marked pyramids and of integer empty tetrahedra are
coincide only for $r=1,2,3,4,5,6,8,10,12,24$.
\end{remark}

\subsection{Proof of Theorem~A for the case of polygonal marked
pyramids}

In this subsection we study the case of marked pyramids with
polygonal bases (containing more than three angles distinct from
the straight angle). In the next subsection we will study
triangular marked pyramids.

\subsubsection{Integer parallelograms contained in some integer
polyhedron}

First of all we show that the integer convex polygons under
consideration contain an integer parallelogram.

\begin{proposition}\label{st1}
Let some four vertices of a convex polygon be integer points.
Then this polygon contains some integer parallelogram that is
integer-affine equivalent either to the standard unit
parallelogram shown on Figure~\ref{2xPar}a$)$, or to the
parallelogram with vertices $(1,0)$, $(0,1)$, $(-1,0)$, and
$(0,-1)$ shown on Figure~\ref{2xPar}b$)$:
\begin{figure}[h]
$$\epsfbox{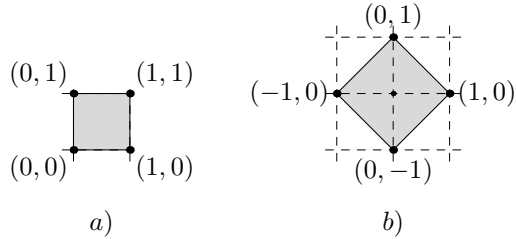}$$
\caption{Any integer polygon contains an integer parallelogram that is integer-affine
equivalent to one of this two parallelograms.}\label{2xPar}
\end{figure}
\end{proposition}

\begin{proof}
Suppose that an integer polygon contains four integer vertices.
Consider the quadrangle generated by these four vertices and denote it by $KLMN$.
Let us prove that the quadrangle contains some integer parallelogram.

Consider the parallelogram $P(KLN)$ and denote it by $KLM'N$.
The vertex $M$ can be in any of the four octants with respect to the lines
containing $M'N$ and $M'L$.
For any of these four cases, we explicitly
construct an integer parallelogram contained in the quadrangle on Figure~\ref{4-gons}
(we draw the quadrangle $KLMN$ with thick line, the corresponding parallelogram is
shaded).
\begin{figure}
$$\epsfbox{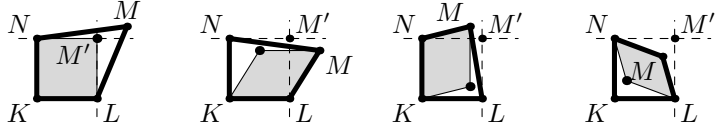}$$
\caption{The possible positions of the convex quadrangle $KLMN$
with respect to $P(KNL)$ (i.e. the quadrangle
$KLM'N$).}\label{4-gons}
\end{figure}

Further we use the following statement.
Let some point of an integer parallelogram be integer.
Consider the point which is centrally-symmetric about the intersection point
of the diagonals of this parallelogram. This point is also
in the parallelogram and is integer.

Denote the integer parallelogram in the polygon by $ABCD$.

{\bf 1. Integer empty parallelogram.} Suppose $ABCD$ is empty.
Then it generates the whole integer lattice and hence is integer-affine
equivalent to the standard one.

{\bf 2. Integer parallelogram with the only one integer point inside.}
Suppose $ABCD$ contains only one integer point $O$ in its interior.
Then this point coincides with the centrally-symmetric point
about the intersection point of the diagonals of this parallelogram.
And hence it coincides with the intersection point of the diagonals.
Therefore the integer triangle $OAB$ is empty. Hence it is
integer-affine equivalent to the standard unit triangle.
Now the integer-affine type of $ABCD$ is uniquely defined and is just
the integer-affine type of the parallelogram with vertices
$(1,0)$, $(0,1)$, $(-1,0)$, and $(0,-1)$.

{\bf 3. Remaining cases.} Let the parallelogram $ABCD$ contain
more than one integer point except the vertices.
Then there exists a points among these points such that it is different
from the intersection point of the diagonals of this parallelogram.
We denote it by$~O$.
Denote the centrally-symmetric point about the intersection point
of the diagonals of this parallelogram by~$O'$.
Without loss of generality, we suppose that~$OO'$ is not a subset of $AC$
(otherwise $OO'$ is not a subset of~$BD$).
Therefore $AOCO'$ (or $AO'CO$) is an integer parallelogram contained in $ABCD$.
The number of integer points of $AOCO'$ is smaller than
the number of integer points of $ABCD$ at least by two.
Since the initial parallelogram contains only a finite number of integer points,
we iteratively come to one of the cases of item {\bf 1.} or {\bf 2.}

Therefore any convex polygon with four integer vertices contains a parallelogram
integer-affine equivalent to one of the parallelograms of Proposition~\ref{st1}.
\end{proof}

\subsubsection{The case of an empty marked pyramid with empty
parallelogram as base}

\begin{proposition}\label{lt1}
Let an empty integer parallelogram be a base of some marked pyramid.
If this pyramid is empty, then it is single-story.
\end{proposition}

\begin{proof}
We prove this proposition by contradiction.
Let $A'ABCD$ be an empty marked pyramid with marked vertex $A'$ and an empty
parallelogram $ABCD$ as its base.
Suppose that the integer distance from the point $A'$ to the plane containing
$ABCD$ equals $r>1$.
Consider the parallelepiped $P(AA'BC)$ and the
integer-distance coordinates corresponding to it
(we denote such coordinates in the following way: $(x,y,z)$).
By Proposition~\ref{m_st}, the coordinates of $A'$, $B$, and $C$
equal to $(r,0,0)$, $(0,r,0)$, and $(0,0,r)$ respectively.
Note that coordinates of lattice nodes (of the integer lattice in the old coordinates)
are integers.

Let us find the lattice node of the parallelepiped
at the unit integer distance to the plane containing $ABC$,
i.e. the lattice node with coordinates $(1,y,z)$, where $0\le y \le r$,
$0 \le z \le r$.
On one hand, it does not contain in the marked pyramid $A'ABCD$,
and hence $y+1>r$ or $z+1>r$.
On the other hand, by Corollary~\ref{l1}, the two-dimensional faces of $P(AA'BC)$
do not contain integer points different from vertices since $AA'BC$ is empty.
Therefore $y$ and $z$ are not equal to $r$.
Hence there are no lattice nodes in the plane containing $ABC$.
We come to the contradiction with Lemma~\ref{l1}.
\end{proof}

\subsubsection{The case of a completely empty marked pyramid whose
base is an integer parallelogram containing a unique integer
point in its interior}

\begin{lemma}\label{mn1}
Consider an integer marked pyramid with vertex $O$ and parallelogram
$ABCD$ as base. Let $ABCD$ be integer-affine equivalent to the parallelogram with
vertices $(1,0)$, $(0,1)$, $(-1,0)$, and $(0,-1)$ $($see Figure~\ref{mnogoug.3}$)$.
If the marked pyramid $OABCD$ is completely empty and multistory, then
it is two-story. The integer-affine type of such pyramid coincides
with the integer-affine type of the pyramid with vertex $(0,0,0)$
and base $(2,-1,0)$, $(2,-2,1)$, $(2,-1,2)$, $(2,0,1)$.
\begin{figure}[h]
$$\epsfbox{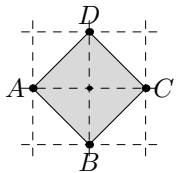}$$
\caption{Quadrangle with vertices $(1,0)$, $(0,1)$, $(-1,0)$, and
$(0,-1)$.}\label{mnogoug.3}
\end{figure}
\end{lemma}

\begin{proof}
Let the integer base $ABCD$ of the completely empty $r$-story
integer marked pyramid $OABCD$ ($r\ge 2$) be integer-affine equivalent to the
parallelogram with vertices $(1,0)$, $(0,1)$, $(-1,0)$, and $(0,-1)$.

Consider the parallelepiped $P(AOBC)$ and the
integer-distance coordinates corresponding to it
(we denote these coordinates as $(x,y,z)$).
By Proposition~\ref{m_st}, the coordinates of $O$, $B$, $C$, and $D$
equal $(r,0,0)$, $(0,2r,0)$, $(0,0,2r)$, and $(0,2r,2r)$ respectively.

Let us consider the parallelogram of intersection of $P(AOBC)$
with
 the plane $x=1$.
Now we will find all lattice nodes in this parallelogram.
By Proposition~\ref{m_st}, there are exactly two lattice nodes in the parallelogram
of intersection.
Let us describe all possible positions of these nodes in the intersection of
$P(AOBC)$ and the plane $x=1$.
First, there are no lattice nodes in the intersection of
the marked pyramid $AOBCD$ and the plane $x=1$,
i.e. in the closed parallelogram with vertices
$(1,0,0)$, $(1,0,2r-2)$, $(1,2r-2, 2r-2)$, and $(1,2r-2,0)$.
Secondly, there are no lattice nodes in all parallelograms
obtained from the given one by applying translations by the vectors
$\lambda (0,2r,0) + \mu (0,r,r)$, where $\lambda$
and $\mu$ are integers.
In Figure~\ref{lemmas.10}, we show some parallelograms
that do not contain any lattice nodes.
These parallelograms are painted shaded.

So, the lattice nodes of the intersection parallelogram
of $P(AOBC)$ with the plane $x=1$ can be only in integer points of
open parallelograms obtained from the parallelogram with vertices
$K(1,r-2,2r-2)$, $L(1,r,2r-2)$, $M(1,r,2r)$, and $N(1,r-2,2r)$
by the symmetry with respect to the plane $y=z$ and translations by the vectors
$\lambda (0,2r,0) + \mu (0,r,r)$, where $\lambda$ and $\mu$ are some integers.
The parallelogram $KLMN$ contains exactly one integer point
$(1,r-1,2r-1)$, see Figure~\ref{lemmas.10}.

\begin{figure}[ht]
$$\epsfbox{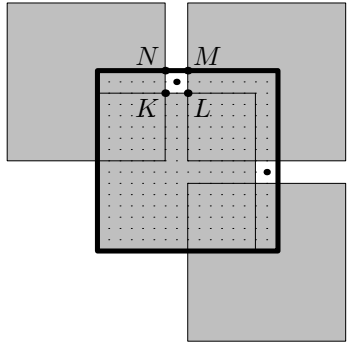}$$
\caption{The intersection of the parallelepiped $P(AOBC)$ and the plane $x=1$.}\label{lemmas.10}
\end{figure}

Suppose that this point is a lattice node.
Since the intersection parallelogram contains exactly two lattice nodes,
the point symmetric to the point $(1,r-1,2r-1)$ with respect to the plane $y=z$  is also a lattice node
(there are no other integer points in the intersection parallelogram).
Therefore $(2,2r-2, 4r-2)$ is a lattice node.
Hence $(2,2r-2,2r-2)$ is a lattice node,
and hence $(2,r-2,r-2)$ is also a lattice node.
However, for $r\ge 3$ the point $(2,r-2,r-2)$ is contained in the closed parallelogram
of intersection of $P(AOBC)$ with the plane $x=2$.
The vertices of this parallelogram are the following:
$(2,0,0)$, $(1,0,2r-4)$, $(1,2r-4, 2r-4)$, and $(1,2r-4,0)$.
Thus there are no pyramids satisfying all the conditions of Lemma~\ref{mn1} for $r\ge 3$.

Now consider the case $r=2$.
The integer points $A$, $B$, $C$, and $(1,1,3)$ define the integer lattice in a unique way.
This implies that all marked pyramids satisfying all the conditions of Lemma~\ref{mn1}
are of the same integer-affine type, and it coincides
with the integer-affine type of the marked pyramid with vertex $(0,0,0)$
and base $(2,-1,0)$, $(2,-2,1)$,
$(2,-1,2)$, $(2,0,1)$ (in the old coordinates).
\end{proof}

\subsubsection{General case}

Now we study the general case of integer completely empty marked
pyramids with convex polygonal bases.

\begin{lemma}\label{mn2}
Consider an integer marked pyramid with vertex $O$ and convex polygonal base $M$.
If this marked pyramid is completely empty and multistory, then
it is two-story. The base of the marked pyramid is integer-affine
equivalent to the quadrangle $(b,0)$, $(0,1)$, $(-a,0)$, $(0,-1)$ where $b\ge a \ge 1$
$($see the case $a=2$, $b=3$ on Figure~\ref{mnogoug.4}$)$.
\begin{figure}[ht]
$$\epsfbox{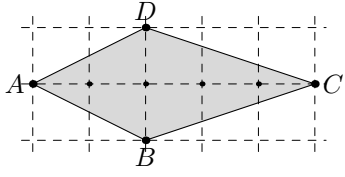}$$
\caption{Quadrangle with vertices $(b,0)$, $(0,1)$, $(-a,0)$, and
$(0,-1)$ where $b\ge a \ge 1$.}\label{mnogoug.4}
\end{figure}
The integer-affine type of the pyramid is uniquely determined
by the integers $a$ and $b$ $($for $b \ge a \ge 1$$)$
and coincides with the integer-affine type of the marked pyramid $M_{a,b}$.
Two marked pyramids $M_{a,b}$ and $M_{a',b'}$ $($$b \ge a \ge 1$, $b' \ge a' \ge 1$$)$
are integer-affine equivalent iff $a=a'$ and $b=b'$.
\end{lemma}

\begin{proof}
Under the assumptions of the lemma, the integer distance from the two-dimensional plane
containing the parallelogram $M$ to the vertex $O$ is greater than one.
From Proposition~\ref{st1} it follows that the parallelogram $M$ contains either
some empty parallelogram or a parallelogram with exactly one integer point in its interior
(and different from the vertices).
By Proposition~\ref{lt1}, the case of an empty parallelogram is eliminated.
Consider an empty parallelogram~$P$ or a parallelogram with exactly one integer point inside.

Choose coordinates on the plane containing the base $M$ so that
the vertices of $P$ have the following coordinates:
$(1,0)$, $(0,1)$, $(-1,0)$, and $(0,-1)$.
Note that all the coordinates of a point of this plane are integers iff
this point is integer (with respect to the old system of coordinates).

Let an integer point with coordinates $(x,y)$ for some $x,y > 0$ be in the base $M$.
Since $M$ is convex, the point $(1,1)$ is also in $M$.
This implies that the empty integer parallelogram with vertices
$(0,0)$, $(1,0)$, $(1,1)$, $(0,1)$ is contained in $M$.
Therefore by, Proposition~\ref{lt1}, the distance from the vertex of the pyramid
to the two-dimensional plane containing
the polygon $M$ equals one.

The cases $x>0$, $y<0$; $x<0$, $y>0$; and $x,y<0$ are similar.

Let the integer points with coordinates $(x,0)$ and $(0,y)$, where $|x|>1$ and $|y|>1$,
be in the base $M$.
Then $M$ contains one of the points: $(1,1)$, $(1,-1)$, $(-1,1)$, or $(-1,-1)$.
And for the same reason, the distance from the vertex of the pyramid
to the two-dimensional plane containing $M$ equals one.

Without loss of generality we suppose that $M$ does not contain points with coordinates
$(0,y)$ for $|y|>1$.
Then $M$ is integer-affine equivalent to the quadrangle with vertices
$(b,0)$, $(0,1)$, $(-a,0)$, $(0,-1)$, where $b\ge a \ge 1$.

Since the polygon $M$ contains the parallelogram $P$,
by Lemma~\ref{mn1} the integer distance from the vertex $O$ of the marked pyramid
to the two-dimensional plane containing the base $M$ equals two.
The parallelogram $P$ is uniquely defined by the quadrangle with vertices
$(b,0)$, $(0,1)$, $(-a,0)$, $(0,-1)$, where $b\ge a \ge 1$
(such a quadrangle contains the unique integer parallelogram
with exactly one integer point different to the vertices).
Therefore, by Lemma~\ref{mn1}, the marked pyramid is integer-affine equivalent
to the marked pyramid with vertex $(0,0,0)$ and base
$(2,-1,0)$, $(2,-a-1,1)$, $(2,-1,2)$, $(2,b-1,1)$.

The point of intersections of the diagonals of the base
quadrangle divides the diagonals into four segments with integer
lengths $1$, $1$, $a$ and $b$. Therefore the (unordered) pair of
integers $[a,b]$ is an integer-affine invariant of the marked
pyramids.
\end{proof}

\subsection{Proof of Theorem~A for the case of triangular marked pyramids}

We continue the proof by exhausting some special cases.
Throughout this subsection we denote by $OABC$ a triangular marked
pyramid with vertex $O$ and base $ABC$.

\subsubsection{Case 1: the base contains an integer polygon}

Suppose that the triangle $ABC$ contains such two integer points
$D$ and $E$, that the line $DE$ intersects the edges of the
triangle $ABC$ and does not contain any vertex of the triangle.
Without loss of generality we suppose that the open ray $DE$ with
vertex at $D$ intersects $AB$, and the open ray $ED$ with vertex
at $E$ intersects $BC$. Hence the triangle $ABC$ contains some
integer convex quadrangle $AEDC$. By Proposition~\ref{st1}, the
triangle $ABC$ contains either an integer empty parallelogram or
a parallelogram integer-affine equivalent to the parallelogram
with vertices $(1,0)$, $(0,1)$, $(-1,0)$, and $(0,-1)$.

If the triangle $ABC$ contains an integer empty parallelogram, then by Proposition~\ref{lt1}
the marked pyramid $OABC$ is single-story.

Suppose that the triangle $ABC$ does not contain an integer empty parallelogram
and contains a parallelogram integer-affine equivalent to the parallelogram with vertices
$(1,0)$, $(0,1)$, $(-1,0)$, and $(0,-1)$.
Consider coordinates on the plane containing the base such that
the vertices of the above-mentioned parallelogram have the following coordinates:
$(1,0)$, $(0,1)$, $(-1,0)$, and $(0,-1)$.
If the points $(1,1)$, $(1,-1)$, $(-1,1)$, and $(-1,-1)$ are not contained
in $ABC$, then the marked pyramid is no longer triangular.
Therefore any marked pyramid of Case~1 contains some empty parallelogram,
and, by Proposition~\ref{lt1}, is single-story.

\subsubsection{Case 2: the integer points of the base different
from the vertices are not contained in one line}

\begin{figure}[h]
$$\epsfbox{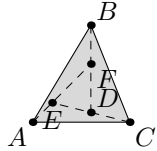}$$
\caption{The affine type of triangles of Case~2.}\label{theorem1.4}
\end{figure}

Now suppose, that there are two integer points $G$ and $H$ such
that the line $GH$ intersects the edges of the triangle $ABC$ and
does not contain any vertex of the triangle. Here we consider the
case of integer points of the base different from the vertices
and not contained in one line. The only possible affine type is
shown on Figure~\ref{theorem1.4}.

Let us find all possible integer-affine types of such a triangle.
Since the triangle $FED$ (see Fig.~\ref{theorem1.4}) is empty,
it is integer-affine equivalent to the triangle $(1,0)$, $(0,0)$, and $(0,1)$.
The points $A$, $B$, and $C$ correspond to $(-1,0)$, $(2,1)$, and
$(0,-2)$ respectively.
Hence the integer-affine type is determined in a unique way.

\begin{lemma}
Consider an integer multistory marked pyramid with vertex $O$ and triangular base $ABC$.
Let the triangle $ABC$ be integer-affine equivalent to the
triangle with vertices $(-2,1)$, $(-1,-1)$, and $(1,2)$,
shown on Figure~\ref{theorem1.5}.
Then the marked pyramid~$OABC$ is two-story and integer-affine equivalent to
the marked pyramid $V$ of List ``M-W''.
\begin{figure}[ht]
$$\epsfbox{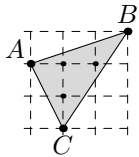}$$
\caption{The triangles with the following vertices $(-2,1)$,
$(-1,-1)$, and $(1,2)$.}\label{theorem1.5}
\end{figure}
\end{lemma}

\begin{proof}
Let the base of an $r$-story ($r \ge 2$) completely empty marked pyramid
$OABC$ be integer-affine equivalent to the
triangle with vertices $(-2,1)$, $(-1,-1)$, and $(1,2)$.

Consider the parallelepiped $P(AOBC)$ and the
integer-distance coordinates corresponding to it
(we denote these coordinates by: $(x,y,z)$).
By Proposition~\ref{m_st}, the coordinates of the vertices $O$, $B$, and $C$
are $(r,0,0)$, $(0,7r,0)$, and $(0,0,7r)$ respectively.

Let us consider the intersection parallelogram of $P(AOBC)$ with
 the plane $x=1$.
Now we will find all lattice nodes in this parallelogram. By
Proposition~\ref{m_st}, there are exactly seven lattice nodes in
the parallelogram of intersection. Let us describe all possible
positions of these nodes in the intersection of $P(AOBC)$ with
the plane $x=1$. First, there are no lattice nodes in the
intersection of the marked pyramid $AOBC$ with the plane $x=1$,
i.e. in the closed triangle with vertices $(1,0,0)$,
$(1,0,7r{-}7)$, and $(1,7r{-}7,0)$. Secondly, there are no
lattice nodes in all triangles obtained from the given one by
applying translations by vectors $\lambda (0,r,2r) + \mu
(0,4r,r)$ for all integers $\lambda$ and $\mu$. In
Figure~\ref{lemmas.1}, ($r\ge 4$) and Figure~\ref{lemmas.2}
($r=2,3$) we show some triangles that do not contain any lattice
nodes. These triangles are shaded.

So the lattice nodes of the intersection parallelogram
of $P(AOBC)$ with the plane $x=1$ can be only at integer points of
open triangles obtained from two triangles by translations by the vectors
$\lambda (0,r,2r) + \mu (0,4r,r)$ for all integers $\lambda$ and $\mu$.
The vertices of the first triangle are
$K(1,3r,4r{-}7)$, $L(1,3r,2r)$, and $M(1,5r{-}7,2r)$.
Here the points $(1,0,0)$ and $L$ should be in different
half-planes with respect to the plane $LM$.
This condition is satisfied only if $2r>4r-7$, i.e. $r< 7/2$.
The vertices of the second triangle are
$P(1,4r-7,3r)$, $Q(1,r,3r)$, and $R(1,r,6r-7)$.
And again the points $(1,0,0)$ and $Q$ should be in different
half-planes with respect to the plane $PR$.
This condition is satisfied only if $(4r-7<r)$, i.e. $r< 7/3$.

So for $r>3$ all points of the intersection parallelogram
of $P(AOBC)$ with the plane $x=1$ are covered, see Figure~\ref{lemmas.1}.
\begin{figure}[ht]
$$\epsfbox{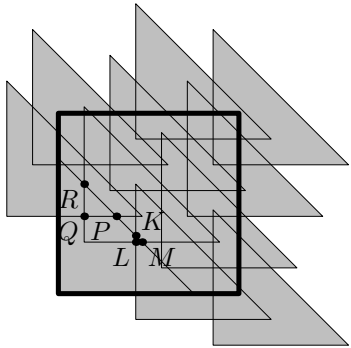}$$
\caption{The intersection of the parallelepiped $P(AOBC)$ with the plane $x=1$ (for $r>3$).}\label{lemmas.1}
\end{figure}
If $r=2$, then the triangle $KLM$ contains only one integer point
with coordinates $(1,5,3)$, see Figure~\ref{lemmas.2}a). If
$r=3$, then the triangle $KLM$ does not contain any integer
point, see Figure~\ref{lemmas.2}b).
\begin{figure}[ht]
$$
\epsfbox{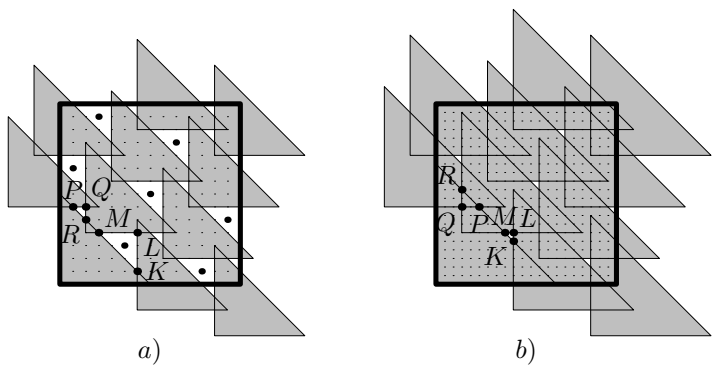}
$$
\caption{The intersection of the parallelepiped $P(AOBC)$ with
the plane $x=1$: a) $r=2$; b) $r=3$.}\label{lemmas.2}
\end{figure}

Since the intersection parallelogram of the plane $x=1$ with the open parallelepiped
must contain seven lattice nodes, the only possible case is the case $r=2$.
There are exactly seven integer points in the complement of the union of the described triangles
in the parallelogram. Hence all these points are lattice nodes.
Therefore, the marked pyramid~$OABC$ is two-story and integer-affine equivalent to
the marked pyramid with vertex $(0,0,0)$ and base $(2,-2,1)$, $(2,-1,-1)$, $(2,1,2)$
(i.e. to the pyramid $V$ of List ``M-W'').
\end{proof}

It remains to study the cases of triangular pyramids with the following property.
All integer points of the base of such a pyramid different from the vertices of the pyramid are
contained in some straight line passing through one of the vertices of the base triangle.

\subsubsection{Case 3: all integer points of the base different
from vertices are contained in a straight line --- {\bf I}}

Suppose that all integer points of the triangle $ABC$ are
contained in a ray with vertex at $A$. Let the number of points
be equal to $c$ ($c\ge 1$), and also suppose all these points are
inner. Denote the inner points by $D_1, \ldots, D_c$, starting
from the point closest to $A$ and increasing the indexing in the
direction from $A$. It turns out that for any positive integer
$c$ there exist exactly one integer-affine type of such pyramid.

\begin{figure}[h]
$$\epsfbox{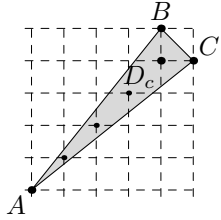}$$
\caption{The triangle with vertices $(1,0)$, $(0,1)$, and $(-c,-c)$.}\label{theorem1.6}
\end{figure}

Since the triangle $BD_cC$ is empty there exists an integer-affine transformation
that maps the triangle to any other empty triangle.
Let us transform the triangle $BD_cC$ to the triangle $\tilde B\tilde D_c\tilde C$
with vertices $(0,1)$, $(0,0)$, and $(1,0)$ respectively.
Now we determine the image of $A$.
Since the point $\tilde D_{c}(0,0)$ is an integer point of the triangle,
the point $\tilde A$ is in the third orthant ($x<0$, $y<0$).
Since $(-1,0)$ is not in the triangle,
the point $\tilde A$ is in the half-plane defined by $y<x+1$.
Since $(0,-1)$ is not in the triangle,
the point $\tilde A$ is in the half-plane defined by $y>x-1$.
Since $\tilde A$ is integer, its coordinates are $(-t,-t)$ for some positive integer $t$.
Since there are exactly $c$ inner integer points in the triangle $\tilde B\tilde D_c\tilde C$,
we obtain $t=c$.
Therefore the triangle $\tilde A \tilde B \tilde C$ is integer-affine equivalent to the triangle
with vertices $(1,0)$, $(0,1)$, and $(-c,-c)$
(the case of $c=4$ is shown on Figure~\ref{theorem1.6}).

First we study the case $c=1$.

\begin{lemma}\label{l3}
Consider an integer multistory marked pyramid with vertex $O$ and triangular base $ABC$.
Let the triangle $ABC$ be integer-affine equivalent to the
triangle with vertices $(-1,-1)$, $(0,1)$, and $(1,0)$
shown on Figure~\ref{theorem1.7}.
Then the marked pyramid~$OABC$ is three-story and integer-affine equivalent to
the marked pyramid $W$ of List ``M-W''.
\begin{figure}[h]
$$\epsfbox{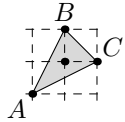}$$
\caption{The triangle with vertices  $(-1,-1)$, $(0,1)$, and $(1,0)$.}\label{theorem1.7}
\end{figure}
\end{lemma}

\begin{proof}
Suppose that the base of $r$-story ($r {\ge} 2$) completely empty
marked pyramid $OABC$ be integer-affine equivalent to the
triangle with the following vertices $(-1,-1)$, $(0,1)$, and
$(1,0)$.

Consider the parallelepiped $P(AOBC)$ and the
integer-distance coordinates corresponding to it
(we denote such coordinates as $(x,y,z)$).
By Proposition~\ref{m_st}, the coordinates of $O$, $B$, and $C$
equal $(r,0,0)$, $(0,3r,0)$, and $(0,0,3r)$ respectively.

Let us consider the parallelogram at intersection of $P(AOBC)$ and
 the plane $x=1$.
Now we will find all lattice nodes in this parallelogram.
By Proposition~\ref{m_st}, there are exactly three lattice nodes in the parallelogram
at intersection.
Let us describe all possible positions of these nodes in the intersection
of $P(AOBC)$ with the plane $x=1$.
First there are no lattice nodes in the intersection of
the marked pyramid $AOBC$ with the plane $x=1$,
i.e. in the closed triangle with vertices
$(1,0,0)$, $(1,0,3r{-}3)$, and $(1,3r{-}3,0)$.
Secondly, there are no lattice nodes in all triangles
obtained from the given one by applying translations by vectors
$\lambda (0,3r,0) + \mu (0,r,r)$ for all integers $\lambda$ and $\mu$.
In Figure~\ref{lemmas.3},
we show some triangles that do not contain any lattice nodes.
These triangles are shaded.

So the lattice nodes of the intersection parallelogram
of $P(AOBC)$ with the plane $x=1$ can be only at integer points of
open triangle obtained from the triangle $K(1,3r,r{-}3)$, $L(1,3r,r)$, $M(1,3r{-}3,r)$
by translations by vectors
$\lambda (0,3r,0) + \mu (0,r,r)$ for all integers $\lambda$ and $\mu$.
Only one point with integer coefficients $(1,3r{-}1,r{-}1)$ is in the triangle $KLM$,
see Figure~\ref{lemmas.3}.

\begin{figure}[ht]
$$\epsfbox{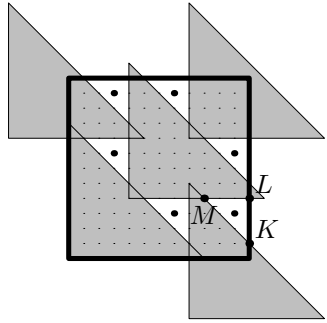}$$
\caption{The intersection of the parallelepiped $P(AOBC)$ with the plane $x=1$.}\label{lemmas.3}
\end{figure}

Shaded triangles covers almost all integer points of the
intersection parallelogram of $P(AOBC)$ with the plane $x=1$.
Only two three-tuples of integer points are still uncovered:\\
{\bf 1)} $(1,3r{-}1,r{-}1)$, $(1,r{-}1,2r{-}1)$, $(1,2r{-}1,3r{-}1)$;\\
{\bf 2)} $(1,r{-}1,3r{-}1)$, $(1,2r{-}1,r{-}1)$, $(1,3r{-}1,2r{-}1)$.\\
So the lattice nodes are either the points of the first
three-tuples or the points of the second one.

Suppose $(1,3r{-}1,r{-}1)$ is a lattice node. (If no, then the
point $(1,r{-}1,3r{-}1)$ is a lattice node. Since the
transformation that maps $(x,y,z)$ to $(x,z,y)$ is integer-affine
and it preserves the parallelepiped $P(AOBC)$ and the marked
pyramid $OABC$, this case is similar.) Then the point
$(r,(3r{-}1)r, (r{-}1)r)$ is a lattice node. Hence
$(3r-1)r-(r-1)r$ is divisible by three, and hence $2r^2$ is also
divisible by three. Therefore $r$ is divisible by three.

Suppose $r=3$, then the marked pyramid exists
and is integer-affine equivalent to $W$.

Let us study the case of $r=3k$, for $k\ge 2$.
Consider the parallelogram at intersection of $P(AOBC)$ and
the plane $x=3$.
Now we will find all lattice nodes in this parallelogram.
By Proposition~\ref{m_st}, there are exactly three lattice nodes in the parallelogram
of intersection.
Let us describe all possible positions of these nodes in the intersection of
of $P(AOBC)$ with the plane $x=3$.
First, there are no lattice nodes in the intersection of
the marked pyramid $AOBC$ with the plane $x=3$,
i.e. in the closed triangle with vertices
$(3,0,0)$, $(3,3r-9,0)$, and $(3,3r-9,0)$.
Secondly, there are no lattice nodes in all triangles
obtained from the given one by applying translations by vectors
$\lambda (0,3r,0) + \mu (0,r,r)$ for all integers $\lambda$ and $\mu$.
This includes the triangle with vertices $P(3,2r,2r)$,
$Q(3,5r-9,2r)$, and $R(3,2r,5r-9)$ shown on Figure~\ref{lemmas.4}
\begin{figure}[ht]
$$\epsfbox{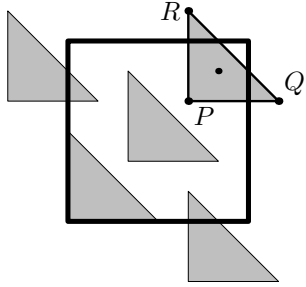}$$
\caption{The intersection of the parallelepiped $P(AOBC)$ with the plane $x=3$.}\label{lemmas.4}
\end{figure}
Since $(1,3r-1,r-1)$ is a lattice node,
the point $(3,9r-3,3r-3)$ is a lattice node.
Thus $(3,3r-3,3r-3)$ is a lattice node.
However, this point is in $KLM$ (for $r>1$) and hence $(1,3r-1,r-1)$ is not a lattice node.
We come to the contradiction, the case of $r=3k$ for $k\ge 2$ is empty.
\end{proof}

\begin{lemma}
Consider an integer multistory marked pyramid with vertex $O$ and triangular base $ABC$.
Let the triangle $ABC$ be integer-affine equivalent to the
triangle with vertices $(-c,-c)$, $(0,-1)$ , and $(-1,0)$, for $c\ge 2$.
Then the marked pyramid~$OABC$ is not completely empty.
\end{lemma}

\begin{proof}
We prove by reductio ad absurdum.
Let the base of $r$-story ($r \ge 2$) completely empty marked pyramid
$OABC$ be integer-affine equivalent to the
triangle with vertices $(-c,-c)$, $(0,-1)$, and $(-1,0)$, for $c\ge 2$.
Since the triangle with vertices
$(-c,-c)$, $(1,0)$, and $(0,1)$ contains the triangle with vertices
$(-1,-1)$, $(1,0)$, and $(0,1)$, the marked pyramid $OABC$ contains
a marked subpyramid integer-affine equivalent to the pyramid
of Lemma~\ref{l3}.
(By {\it marked subpyramid $P$} of some marked pyramid $Q$ we call such convex
pyramid $P$ that the vertices of $P$ and $Q$ coincides and
the base of $Q$ contains the base of $P$.)
Therefore by Lemma~\ref{l3} we have $r=3$.

Let us show that $r\ne 3$.
Suppose $r=3$.
Since $c\ge 2$, the marked pyramid $OABC$ contains some
marked subpyramid $OA'BC$ with base $A'BC$ integer-affine equivalent to
the triangle with vertices $(-2,-2)$, $(1,0)$, and $(0,1)$.
We show now that $OA'BC$ is not completely empty.

Consider the parallelepiped $P(A'OBC)$ and the
integer-distance coordinates corresponding to it
(we denote such coordinates as $(x,y,z)$).
By Proposition~\ref{m_st}, the coordinates of $O$, $B$, and $C$
equal $(3,0,0)$, $(0,15,0)$, and $(0,0,15)$ respectively.

Let us consider the parallelogram at intersection of $P(A'OBC)$
and the plane $x=1$. Now we will find all lattice nodes in this
parallelogram. First, there are no lattice nodes in the
intersection of the marked pyramid $A'OBC$ with the plane $x=1$,
i.e. in the closed triangle with vertices $(1,0,0)$, $(1,0,12)$,
and $(1,12,0)$. Secondly, there are no lattice nodes in all
triangles obtained from the given one by applying translations by
vectors $\lambda (0,15,0) + \mu (0,3,3)$ for all integers
$\lambda$ and $\mu$. This triangles contains all integer points
of the intersection of $P(A'OBC)$ with the plane $x=1$, see
Figure~\ref{lemmas.5}.

\begin{figure}[h]
$$\epsfbox{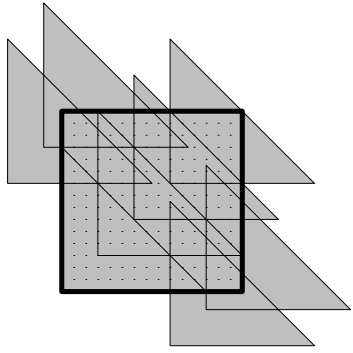}$$
\caption{The intersection of the parallelepiped $P(A'OBC)$ with the plane $x=1$.}\label{lemmas.5}
\end{figure}

So, the marked pyramid $OA'BC$ is not completely empty.
Hence the marked pyramid $OABC$ is not completely empty.
Thus $r\ne 3$.

Therefore, for any $r \ge 2$, the base of any $r$-story completely empty
pyramid $OABC$ is not integer-affine equivalent to the
triangle with vertices $(-c,-c)$, $(0,-1)$, and $(-1,0)$, for $c\ge 2$.
We come to the contradiction.
\end{proof}

\subsubsection{Case 4: all integer points of the base different
from vertices are contained in a straight line --- {\bf II}}

Suppose that all integer points of the triangle $ABC$ are contained in
the ray with vertex $A$. Let the number of points be equal to $b$ ($b\ge 1$),
and the last point is in the edge $BC$.
Denote these points by $D_1, \ldots, D_c$, starting
from the point closest to $A$ and increasing the indexing in the direction from $A$.
It turns out that for any $b$ there exist exactly one integer-affine type of such
pyramid.

\begin{figure}[ht]
$$\epsfbox{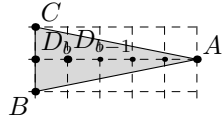}$$
\caption{The triangle with vertices $(0,-1)$, $(b,0)$, and $(0,1)$.}\label{theorem1.8}
\end{figure}

Since the triangle $D_b D_{b-1} B$ is empty there exists an integer-affine transformation
that maps the triangle to any other empty triangle.
We transform the triangle $D_b D_{b-1} B$ to the triangle with vertices
$(0,0)$, $(1,0)$, and $(0,-1)$ respectively.
Then $C$ maps to $(0,1)$, and $A$ maps to $(b,0)$.
Therefore the triangle $ABC$ is integer-affine equivalent to the triangle
with vertices $(0,-1)$, $(b,0)$, and $(0,1)$.
(the case of $b=5$ is shown on Figure~\ref{theorem1.8}).

First we study the case $b=2$.

\begin{lemma}\label{l5}
Consider an integer multistory marked pyramid with vertex $O$ and triangular base $ABC$.
Let the triangle $ABC$ be integer-affine equivalent to the
triangle with vertices $(2,0)$, $(0,-1)$, and $(0,1)$
shown on Figure~\ref{theorem1.9}.
Then the marked pyramid~$OABC$ is two-story and integer-affine equivalent to
the marked pyramid $U_2$ of List ``M-W''.
\begin{figure}[h]
$$\epsfbox{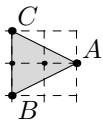}$$
\caption{The triangle with vertices  $(2,0)$, $(0,-1)$, and $(0,1)$.}\label{theorem1.9}
\end{figure}
\end{lemma}

\begin{proof}
Suppose that the base of $r$-story ($r \ge 2$) completely empty
marked pyramid $OABC$ be integer-affine equivalent to the
triangle with vertices $(2,0)$, $(0,-1)$, and $(0,1)$.

Consider the parallelepiped $P(AOBC)$ and the
integer-distance coordinates corresponding to it
(we denote such coordinates as $(x,y,z)$).
By Proposition~\ref{m_st}, the coordinates of $O$, $B$, and $C$
equal $(r,0,0)$, $(0,4r,0)$, and $(0,0,4r)$ respectively.

Now consider the parallelogram at intersection of $P(AOBC)$ and
 the plane $x=1$.
Now we will find all lattice nodes in this parallelogram.
By Proposition~\ref{m_st}, there are exactly three lattice nodes in the parallelogram
at intersection.
Let us describe all possible positions of these nodes in the intersection
of $P(AOBC)$ with the plane $x=1$.
First, there are no lattice nodes in the intersection of
the marked pyramid $AOBC$ with the plane $x=1$,
i.e. in the closed triangle with vertices
$(1,0,0)$, $(1,0,4r-4)$, and $(1,4r-4,0)$.
Secondly, there are no lattice nodes in all triangles
obtained from the given one by applying translations vectors
$\lambda (0,4r,0) + \mu (0,r,r)$ for all integers $\lambda$ and $\mu$.
We show (shaded) triangles that do not contain any lattice nodes on Figure~\ref{lemmas.7}.

So the lattice nodes of the intersection parallelogram
of $P(AOBC)$ with the plane $x=1$ can be only at integer points of
open triangle obtained from the triangle $K(1,4r,2r-3)$, $L(1,4r,2r)$, $M(1,4r-3,2r)$
by translations by vectors
$\lambda (0,4r,0) + \mu (0,r,r)$ for all integers $\lambda$ and $\mu$
and the symmetry about the plane $y=z$.
Only the points with integer coefficients $(1,4r-2,2r-1)$, $(1,4r-1,2r-1)$, and $(1,4r-1,2r-2)$
are in the triangle $KLM$, see Figure~\ref{lemmas.7}.

\begin{figure}[ht]
$$\epsfbox{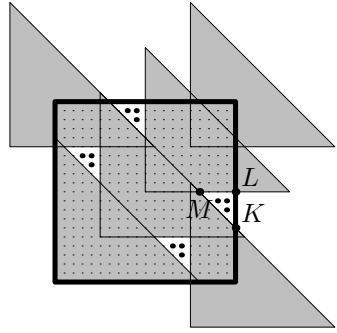}$$
\caption{The intersection of the parallelepiped $P(AOBC)$ with the plane $x=1$.}\label{lemmas.7}
\end{figure}

We prove that one of these points is a lattice node by reductio
ad absurdum. Suppose that the triangle $KLM$ does not contain a
lattice node. Then there are no lattice nodes in all triangles
obtained from $KLM$ by applying translations by vectors of the
form $\lambda (0,4r,0) + \mu (0,r,r)$ for all integers $\lambda$
and $\mu$. Hence the intersection of the parallelepiped $P(AOBC)$
with the plane $x=1$ does not contain integer nodes. We come to
the contradiction. So one of the points $(1,4r-2,2r-1)$,
$(1,4r-1,2r-1)$, and $(1,4r-1,2r-2)$ is a lattice node.

Suppose that $r \ge 3$ and consider the plane $x=2$.
First, there are no lattice nodes in the intersection of
the marked pyramid $AOBC$ with the plane $x=2$,
i.e. in the closed triangle with vertices
$(1,0,0)$, $(1,0,4r-8)$, and $(1,4r-8,0)$.
Secondly, there are no lattice nodes in all triangles
obtained from the given one by applying translations by vectors
$\lambda (0,4r,0) + \mu (0,r,r)$ for all integers $\lambda$ and $\mu$.
In particular, there are no lattice nodes in the triangle
with vertices $P(2,3r,3r)$,
$Q(2,7r-8,3r)$, and $R(2,3r,7r-8)$.
\begin{figure}[ht]
$$\epsfbox{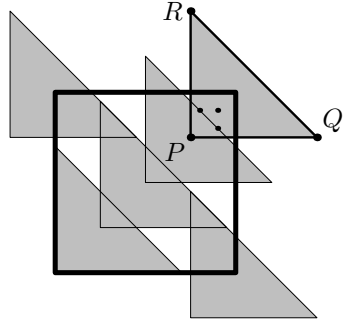}$$
\caption{The intersection of the parallelepiped $P(AOBC)$ with the plane $x=2$.}\label{lemmas.8}
\end{figure}

Suppose that the point
$(1,4r-2,2r-1)$, $(1,4r-1,2r-1)$, or $(1,4r-1,2r-2)$ is a lattice node,
then
$(2,8r-4,4r-2)$, $(2,8r-2,4r-2)$, or $(2,8r-2,4r-4)$ respectively is also a lattice node.
Hence the point
$(2,4r-4,4r-2)$, $(2,4r-2,4r-2)$, or $(2,4r-2,4r-4)$ respectively is a lattice node.
The last three points are contained in the triangle $PQR$
with vertices
$P(2,3r,3r)$, $Q(2,7r-8,3r)$, and $R(2,3r,7r-8)$, for $r>3$ (see Figure~\ref{lemmas.8}),
and hence these points are not lattice nodes.
For $r=3$, the point $(1,11,5)$ is not a lattice node by the same reason.
The points $(1,10,5)$ and $(1,11,4)$ are not lattice nodes,
since the points $(3,30,15)$ and $(3,33,12)$ are not lattice nodes of the plane $x=3$
(all such node coordinates are $(3,4m,4n)$ for some integers $m$ and $n$).
From the above we conclude that $r\le 2$.

Suppose now that $r{=}2$ and consider the points $(1,6,4)$,
$(1,7,3)$, and $(1,7,4)$. The points $(1,6,4)$ and $(1,7,3)$ are
not lattice nodes since the points $(2,12,6)$ and $(2,14,8)$ are
not lattice nodes of the plane $x{=}2$ (all such nodes coordinates
are $(2,4m,4n)$ for some integers $m$ and $n$). The point
$(1,7,4)$ defines a unique-possible integer-affine type of marked
pyramids with such base --- the integer-affine type of the marked
pyramid $U_2$.
\end{proof}

Now we will study the general case ($b\ge 2$).

\begin{lemma}
Consider an integer multistory marked pyramid with vertex $O$ and triangle base $ABC$.
Let the triangle $ABC$ be integer-affine equivalent to the
triangle with vertices $(b,0)$, $(0,-1)$, and $(0,1)$, for $b\ge 2$.
Then the marked pyramid~$OABC$ is two-story and integer-affine equivalent to
the marked pyramid $U_b$ of List ``M-W''.
\end{lemma}

\begin{proof}
Let the base of $r$-story ($r {\ge} 2$) completely em\-pty marked
pyramid $OABC$ be integer-affine equi\-valent to the triangle with
vertices $(b,0)$, $(0,-1)$, and $(0,1)$.

Since the triangle with vertices $(b,0)$, $(0,-1)$, and $(0,1)$
contains the triangle with vertices $(2,0)$, $(0,-1)$, and $(0,1)$,
the marked pyramid $OABC$ contains a marked subpyramid that is
integer-affine equivalent to the marked pyramid of Lemma~\ref{l5}.
Since the subpyramid is completely empty, by Lemma~\ref{l5}
we have that it is two-story.

Suppose now $r{=}2$. Consider the parallelepiped $P(AOBC)$ and the
integer-distance coordinates corresponding to it (let us denote
such coordinates as $(x,y,z)$). By Proposition~\ref{m_st}, the
coordinates of $O$, $B$, and $C$ equal $(2,0,0)$, $(0,4b,0)$, and
$(0,0,4b)$ respectively.

Consider the parallelogram at the intersection of $P(AOBC)$ and
 the plane $x=1$.
Now we will find all lattice nodes in this parallelogram.
By Proposition~\ref{m_st}, there are exactly $2b$ lattice nodes in the parallelogram
at intersection.
Let us describe all possible positions of these nodes in the intersection
of $P(AOBC)$ with the plane $x=1$.
First, there are no lattice nodes in the intersection of
the marked pyramid $AOBC$ with the plane $x=1$,
i.e. in the closed triangle with vertices
$(1,0,0)$, $(1,0,2b)$, and $(1,2b,0)$.
Secondly, there are no lattice nodes in all triangles
obtained from the given one by applying translations by vectors
$\lambda (0,4b,0) + \mu (0,2,2)$ for all integers $\lambda$ and $\mu$.
We show some (shaded) triangles that do not contain any lattice nodes on Figure~\ref{lemmas.9}.

So the lattice nodes of the intersection parallelogram
of $P(AOBC)$ with the plane $x=1$ can be only at integer points of
open triangle obtained from the triangle $K(1,4b,2b-4)$, $L(1,4b,2b)$, $M(1,4b-4,2b)$
by translations by vectors
$\lambda (0,4b,0) + \mu (0,2,2)$ for all integers $\lambda$ and $\mu$
and the symmetry about the plane $y=z$.
Only the points with integer coefficients $(1,4b-2,2b-1)$, $(1,4b-1,2b-1)$, and $(1,4b-1,2b-2)$
are in the triangle $KLM$ (the case $b=3$  is shown on Figure~\ref{lemmas.9}).

\begin{figure}[ht]
$$\epsfbox{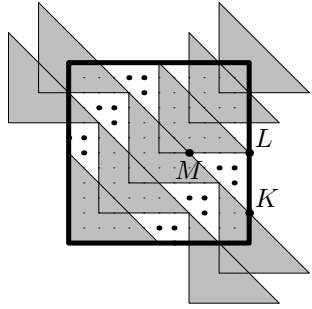}$$
\caption{The intersection of the parallelepiped $P(AOBC)$ with the plane $x=1$.}\label{lemmas.9}
\end{figure}

One of the integer points of this triangle is a lattice node (the
other uncovered parts of the section can be obtained by
translations by vectors $\lambda (0,4b,0) + \mu (0,2,2)$ for all
integers $\lambda$ and $\mu$).

Consider the plane $x=2$.
The point $(2,y,z)$ is a lattice node iff there exist such integers
$m$ and $n$ that $z=2m$, and $y=2m+2bn$.

We show that the point $(1,4b{-}2,2b{-}1)$ is not a lattice node by reductio ad absurdum.
Let this point be a lattice node.
Then the point $(2,8b{-}4,4b{-}2)$ is also a lattice node.
Let us find the such integers $m$ and $n$ that $4b-2=2m$ and $8b-4=2m+2bn$.
Then $m=2b-1$, $n=\frac{2b{-}1}{b}$. For $b \ge 2$, the number $n$ is not integer.
We come to the contradiction. Therefore the point $(1,4b{-}2,2b{-}1)$
is not a lattice node.

By the same reasons the point $(1,4b-1,2b-2)$ is not a lattice node.
The last point of the triangle $(1,4b-1,2b-1)$ determines the pyramid of the
integer-affine type $U_b$.
\end{proof}

\subsubsection{Case 5: integer points of the base different from
vertices are contained in one edge of the base}

\begin{figure}[ht]
$$\epsfbox{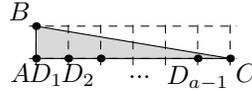}$$
\caption{The triangle with vertices $(0,0)$, $(0,1)$, and $(a,0)$.}\label{theorem1.10}
\end{figure}

It remains to study the case of the last most simple series of triangular
marked pyramids.
Suppose that all integer points of the base $ABC$ different from the vertices
are contained in $AC$, and the integer length of $AC$ is $a-1$, for some $a\ge 2$.
The case of $a=1$ is the case of empty marked pyramid was studied before in
Corollary~\ref{col_pyr2}.
Denote these points by $D_1, \ldots, D_{a-1}$ starting
from the point closest to $A$ and increasing the indexing in the direction to $C$.
(See Figure~\ref{theorem1.10}.)

Consider an integer multistory marked pyramid with vertex $O$ and triangular base $ABC$.
Let the triangle $ABC$ be integer-affine equivalent to the
triangle with vertices $(0,0)$, $(0,1)$, and $(a,0)$, for $a\ge 2$.

Then the marked pyramid~$OABC$ is two-story and integer-affine equivalent to
the marked pyramid $U_b$ of List ``M-W''.

\begin{lemma}\label{zzzzz}
The marked pyramid~$OABC$ is integer-affine equivalent to
the marked pyramid of the following list.

{\bf List ``T'':}

--- $T_{a,1}^0$;

--- $T_{a,r}^\xi$, where
$\xi$ and $r$ are relatively prime and satisfy: $r \ge 2$ and
$0<\xi \le r/2$.

All integer marked pyramids listed in $``T''$ are completely empty and
integer-linear nonequivalent to each other.
\end{lemma}

\begin{proof}
{\bf 1. Preliminary statement.}
Let us show that the marked pyramid $OABC$ is integer-affine equivalent
to the marked pyramid $T_{a,r}^{\xi}$, for some positive integer $\xi \le r/2$.

First of all two single-story marked pyramids with the same $a$
are integer-affine equivalent, since the integer points of the edges of the pyramid
generates all integer lattice.

Let the base of $r$-story ($r \ge 2$) completely empty marked pyramid
$OABC$ be integer-affine equivalent to the
triangle with vertices $(0,0)$, $(0,1)$, and $(a,0)$.
Consider the parallelepiped $P(AOBD_1)$ and the
integer-distance coordinates corresponding to it
(we denote such coordinates as $(x,y,z)$).
By Proposition~\ref{m_st}, the coordinates of $O$, $B$, and $C$
equal $(r,0,0)$, $(0,r,0)$, and $(0,0,r)$ respectively.

By Corollary~\ref{symplex} (since the tetrahedron $AOBD_1$ is empty)
all inner lattice nodes are contained in one of three  diagonal planes:
$x+z=r$, $y+z=r$, or $x+y=r$. Examine all the cases.

Let all inner lattice nodes are contained in the plane $x+z=r$.
By Lemma~\ref{l1} there exist exactly one lattice node $K$ contained in the
plane $x=1$.
So, $K$ is in the intersection of these two planes,
and its coordinates are $(1,\xi,r-1)$, where $0<\xi<r$.
Now we come back to the old coordinates associated with the lattice.
Since the integer distance from $K$ to the two-dimensional plane
containing the face $AD_1B$ equals one,
the tetrahedron $AD_1BK$ can be mapped by some integer-affine transformation
to the tetrahedron with vertices $(0,0,0)$,
$(1,0,0)$, $(0,1,0)$, and $(0,0,1)$.
By such transformation the vertex $O$ maps to $(-\xi,1-r,r)$, and $C$
maps to $(a,0,0)$.
Let us translate the obtained pyramid by the integer vector $(\xi , r-1,r)$.
Finally we get the marked pyramid $T_{a,r}^{\xi}$.
Hence the marked pyramid $OACB$ is integer-affine equivalent to the marked pyramid
$T_{a,r}^\xi$, where $0<\xi< r$.
Consider the integer-affine transformation mapping the points
$O$, $A$, $B$, $C$ to the points $O$, $C$, $B$, $A$ respectively,
then the point $K$ maps to the point $(r-\xi,1-r,r)$.
Chose the smallest one of $\xi$ and $r-\xi$. Obviously,
this number is less then $r/2$.

Let all inner lattice nodes be contained in the plane $y+z=r$
in the integer-distance coordinate system.
By Lemma~\ref{l1} there exists exactly one lattice node $K$ contained in the
plane $x=1$.
So, $K$ is in the intersection of these two planes,
and its coordinates are $(1,\xi,r-\xi)$, where $0<\xi<r$.
The intersection of the marked pyramid $OABC$ with the plane $x=1$
is a triangle with vertices $(1,0,0)$, $(1,ar-a,0)$, and $(1,0,r-1)$.
This triangle contains all integer points $(1,t,r-t)$, for $2\le t\le r$. Hence
$\xi=1$.
Therefore the point $K$ is in the plane $x+z=r$,
so, we are in the position of the previous case.

Let all inner lattice nodes be contained in the plane $x+y=r$
in the integer-distance coordinate system.
By Lemma~\ref{l1} there exist exactly one lattice node $K$ contained in the
plane $z=1$.
So, $K$ is in the intersection of these two planes,
and its coordinates are $(\xi,r-\xi,1)$, where $0<\xi<r$.
The intersection of the marked pyramid $OABC$ with the plane $z=1$
is a triangle with vertices $(0,0,1)$, $(r-1,0,1)$, and $(0,ar-a,1)$.
This triangle contains all integer points $(t,r-t,1)$, for $1\le t\le r-1$. Hence
 $\xi=r-1$.
Therefore the point $K$ is in the plane $x+z=r$,
so, we are in the position of the previous case.

So, the marked pyramid~$OABC$ is integer-affine equivalent to
a marked pyramid $T_{a,r}^{\xi}$, for some positive integer $\xi \le r/2$.

{\bf 2. Completeness of List ``T'' and completely emptiness
of the marked pyramids of ``T''.}
Let us show that the marked pyramids $T_{a,r}^{\xi}$ of the list ``T''
are completely empty.
Denote the vertices of the marked pyramids by $O$, $A$, $B$, $C$,
and the integer points of $AC$ by $D_i$.

Denote also the point $A$ by $D_0$, and the point $C$ by $D_a$.
Note that the marked pyramid $OD_{i}D_{i+1}B$ is integer-affine equivalent
to the marked pyramid $P_r^{\xi}$, for any positive integer $i\le a$,
since the marked pyramid $OD_{i}D_{i+1}B$ can be obtained from the
pyramid $P_r^{\xi}$ by applying the compositions of
the integer-linear transformation defined by the following matrix
$$
\left(
\begin{array}{ccc}
\xi +i+1& \xi +i& -\xi -i\\
r-1& r-1& 2-r\\
-r&-r&r-1\\
\end{array}
\right)
,
$$
and the translation by the integer vector $(-\xi, 1-r,r)$.

By corollary~\ref{col_pyr}, if $\xi$ and $r$ are relatively prime, then
the marked pyramids $OAD_1B$, $OD_1D_2B$,~$\ldots$, $OD_{a-1}CB$
are empty, and hence their union $OABC$ is completely empty.

By the same reasons the marked pyramids $T_{a,r}^{\xi}$ with
relatively prime $\xi$ and $r$ are completely empty.

Therefore List ``T'' is complete, and all integer pyramids of the list
are completely empty.

{\bf 3. Irredundance of List~``T''.}
Now we prove that all marked pyramids $T^\xi_{a,r}$ of List~``T''
are integer-affine nonequivalent to each other.
Obviously, that the marked pyramids with different $a$ are nonequivalent.
Since the integer distance from the marked vertex to the two-dimensional
plane of the marked base is an integer-affine invariant,
the marked pyramids with different $r$ are nonequivalent.

For the case of pyramids with the same integers $a{>}1$ and $r$,
we construct the follo\-wing integer-linear invariant. Consider an
arbitrary marked pyramid $OABC$, where all its integer points are
contained in the edge $AC$. As it was shown before the empty
marked pyramids $OAD_1B$, $OD_1D_2B$,~$\ldots$, $OD_{a-1}CB$ are
integer-affine equivalent to the mar\-ked pyramid $P_r^\xi$ with
$0\le\xi\le r/2$. Since the col\-lection of this marked pyramids
is defined in a unique way and by Corollary~\ref{col_pyr}, the
type of such $P_r^\xi$ is an invariant. This invariant
distinguishes different marked pyramids of List~``T''.
\end{proof}

So, we have studied all possible cases of integer-affine types of
multistory completely empty convex three-dimensional marked
pyramids. It remains to say a few words about the irredundance of
List ``M-W'' of Theorem~A.

\subsubsection{Irredundance of List ``M-W''}

If two marked pyramids have integer-affine nonequivalent bases, then
these pyramids are also integer-affine nonequivalent.
The integer-affine types of the base distinguish almost all marked
pyramids of List ``M-W''.
This does not work only for pyramids $T_{a,r}^{\xi}$ with the same $a$ and $r$,
and distinct $\xi$ from List~``M-W''.
Such pyramids $T_{a,r}^{\xi}$ are integer-affine nonequivalent by Lemma~\ref{zzzzz}
(see List ``T'').

The proof of the main theorem is completed.
\qed

\section{Proof of Theorem~B}\label{pB}

\subsection{Completeness of Lists ``$\alpha_n$'' for $n \ge 2$ of Theorem~B }

Note that polygonal faces of any sail are faces of the boundary of the convex.
Hence all faces are convex.
Consider some marked pyramid with marked vertex at the origin
and some compact two-dimensional face of a sail as base.
By the definition of multidimensional continued fractions it follows
that such pyramid is completely empty.

\begin{lemma}\label{l81}
Two two-dimensional faces are integer-linear equivalent iff
the corresponding completely empty marked pyramids are integer-affine equivalent.
\end{lemma}

\begin{proof}
If two two-dimensional faces are integer-linear equ\-iva\-lent,
then one of them maps to the other with some integer-line\-ar transformation.
The marked pyramid corresponding to the first face maps
to the marked pyramid corresponding to the second face at that.

Suppose now that the corresponding completely empty marked
pyramids are in\-te\-ger-affine equivalent. Then one of them maps
to the other with some integer-affine transformation. Since the
marked vertices of both pyramids are at the origin, the origin is
a fixed point of the transformation. Hence the transformation is
integer-linear. Since the base of the first pyramid maps to the
base of the second, the first face maps to the second also. Hence
these two-dimensional faces are integer-linear equivalent.
\end{proof}

So, for any $r\ge 2$, the following is true. Any integer-linear
type of compact two-dimensional faces contained in the
two-dimensional planes at integer distances equal $r$ from the
origin is uniquely defined by the corresponding integer-affine
type of $r$-story completely empty convex marked pyramids. Hence
by Theorem~A (see List~``M-W'') Lists~``$\alpha_n$'' of theorem~B
are complete if $n>2$.

Now we study the case of two-dimensional continued fractions. By
Theorem~A the list of all triangular faces in List~``$\alpha_2$''
is complete. It remains to show that there are no faces of sails
integer-linear equivalent to the quadrangle with vertices
$(2,-1,0)$, $(2,-a-1,1)$, $(2,-1,2)$, $(2,b-1,1)$, for some $b
\ge a \ge 1$.

Let us prove the following lemma.

\begin{lemma}\label{zapret}
{\bf (On some restrictions for two-di\-men\-sional continued
fractions.)} Let some compact two-dimensional face $F$ of some
two-dimensional continued fraction contains some integer
parallelog\-ram~$P$ integer-affine equivalent to the
parallelogram the with vertices
 $(1,0)$, $(0,1)$, $(-1,0)$, and $(0,-1)$.
Then the integer distance from the origin to the plane containing  the face $F$
equals one.
\end{lemma}

\begin{proof}
Consider such coordinates on the plane containing $F$ that
the vertices of $P$ become $(1,0)$, $(0,1)$, $(-1,0)$, and $(0,-1)$.
Note that the point in this plane is integer iff
its new coordinates are integers.

Suppose that the point $(1,1)$ is in $F$. Then the empty
parallelogram with vertices $(0,0)$, $(1,0)$, $(1,1)$, and
$(0,1)$ is contained in $F$. Therefore, by Proposition~\ref{lt1},
the distance from the plane containing $F$ equals one.

By the same reasons if $(1,-1)$, or $(-1,1)$, or $(-1,-1)$ is in
$F$, then the distance from the plane containing $F$ equals one.

Now we show that $F$ contains one of the listed four points by reduction ad absurdum.
Suppose the points $(1,1)$, $(1,-1)$, $(-1,1)$, and $(-1,-1)$
are in the complement to $F$.
Three planes of two-dimensional continued fraction intersect with the plane
containing $F$ at three lines. The face $F$ is in the interior the triangle $T$
generated by the intersection lines.
The triangle $T$ contains $F$, and the set $T\setminus F$ does not
contain any integer point.
Notice that the point $(1,0)$ is in $F$,  and the points $(1,1)$ and $(1,-1)$
are not in $F$. Note also that the points $(1,0)$, $(1,1)$, and $(1,-1)$
are in one straight line.
Then the open angle with vertex $(0,0)$ and edges passing through
the points $(1,1)$, and $(1,-1)$,
contains some vertex of the triangle $T$, see Figure~\ref{mnogoug.2}.

\begin{figure}[ht]
$$\epsfbox{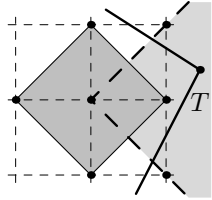}$$
\caption{One of the vertices of $T$ is in shaded
(open) angle.}\label{mnogoug.2}
\end{figure}

The same holds for two adjacent angles and for the opposite angle.
Therefore the triangle $T$ has at least four vertices.
We come to the contradiction.

So, we have studied all the cases. Lemma~\ref{zapret} is proven.
\end{proof}

\begin{corollary}
Any two-dimensional continued fractions does not contain faces
that are integer-linear equivalent to the quadrangle with vertices
$(2,-1,0)$, $(2,-a{-}1,1)$, $(2,-1,2)$, $(2,b{-}1,1)$ for $b \ge a \ge 1$.
\end{corollary}

\begin{proof}
Take the polygon with the following vertices $(2,-1,0)$,
$(2,-a-1,1)$, $(2,-1,2)$, $(2,b-1,1)$. It contains the
parallelogram integer-affine equivalent to the parallelogram with
vertices $(1,0)$, $(0,1)$, $(-1,0)$, and $(0,-1)$. Hence by
Lemma~\ref{zapret} the distance between the origin and the plane
containing such face equals one. But the integer distances from
the origin to the faces of the corollary equal two. Therefore
two-dimensional continued fractions do not contain faces of the
corollary.
\end{proof}

Therefore all Lists ``$\alpha_n$'' of Theorem~B (for
$n$-dimensional continued fractions) are completed for all $n\ge
2$.

\subsection{The completion of proof of Theorem~B}

 First we show that all triangular faces of List ``$\alpha_n$''
are realizable for $n\ge 2$. We prove more general statement for
the triangles, and then generalize it to the case of polygons.

\subsubsection{All triangular faces are realizable}

Consider some completely empty triangular marked pyramid $OABC$
with marked vertex $O$ and base $ABC$. Let $\Sigma_{ABC}(3)$ be
the configuration space of ordered 3-tuples of points of the
plane containing the triangle $ABC$. The configuration space
$\Sigma_{ABC}(3)$ is homeomorphic to $\r ^6$. Consider the
standard topology of this space. Denote by $U_3 \subset
\Sigma_{ABC}(3)$ the set of such 3-tuples $A'B'C'$, that the open
marked pyramid $OA'B'C'$ contains the marked pyramid $OABC$
(except the point $O$) and the set $OA'B'C'\setminus OABC$ does
not contain integer points.

\begin{lemma}\label{last}
The set $U_3$ is open and nonempty.
Any point of the set $U_3$ defines such two-dimensional continued fraction,
that this fraction contains
the triangle $ABC$ as a two-dimensional face.
\end{lemma}

\begin{proof}
First we prove that $U_3$ is open. Consider such integer planes
parallel to the plane $ABC$ that the origin and the pane $ABC$
are in different half-planes (in the three-dimensional space
spanning the points $O$, $A$, $B$, and $C$). The number of such
planes is finite and equals $r-1$, where $r$ is an integer
distance between $O$ and the plane $ABC$. Denote by $\pi_i$, for
$i\le r$, one of the described planes at the integer distance
from $O$ equal~$i$. The marked pyramid intersects $OABC$ with
$\pi_i$ by the triangle, we denote it by $T_i$. The triangle
$T_i$ does not contain integer points for $i<r$. Consider all
open triangular angles centered at $O$, that intersect with
$\pi_i$ by some triangles that contain closed $T_i$ and do not
contain other integer points (different from the integer points
of $T_r$ for the case of $i=r$). Any such angle defines three
points in the plane containing the triangle $ABC$ (i.e. in
$\pi_r$). These points determine six ordered 3-tuples points of
$U_3$. The set of all such triangular angles determines the
nonempty open subset of $U_3$ (denoted by $U_{3,i}$). As it is
easy to see, the set $U_3$ coincides with  the intersection of
the sets $U_{3,i}$ for all positive integers $i\le r$. Therefore
$U_3$ is open.

Secondly we prove that $U_3$ is nonempty.
We denote by $u_0$ the point of $U_3$ corresponding to the ordered 3-tuple
points $A$, $B$, $C$.
On one hand, there exist a neighborhood of $u_0$
containing the points with the following property:
if $A'B'C'$ in the neighborhood, then the set of
integer points of the marked pyramid $OA'B'C'$
is contained in the marked pyramid $OABC$.
On the other hand, any neighborhood of $u_0$ (and also the described one)
contains such point $A''B''C''$ that the closed marked pyramid $OABC$ is contained
in the open triangular angle $OA''B''C''$ except the point $O$.
From these two facts it follows that $U_{3}$ is not empty.

Now consider an arbitrary triangle $A'B'C'$ corresponding to some
point of $U_3$. Chose the planes $OA'B'$, $OA'C'$, and $OB'C'$.
The two-di\-men\-si\-onal continued fraction defined by these
planes contains the triangle $ABC$ as a face.
\end{proof}

\subsubsection{Realizability of polygonal faces}

Let us now generalize Lemma~\ref{last}. Consider some completely
empty convex polygonal  marked pyramid $OA_1\ldots A_n$ with
marked vertex $O$ and base $A_1\ldots A_n$. Let $\Sigma_{A_1
\ldots A_n}(n)$ be the configuration space of ordered $n$-tuples
of points of the plane containing the polygon $A_1\ldots A_n$.
The configuration space $\Sigma_{A_1 \ldots A_n}(n)$ is
homeomorphic to $\r ^{2n}$. Consider the standard topology of
this space. Denote by $U_n \subset \Sigma_{A_1 \ldots A_n}(n)$
the set of such $n$-tuples $A'_1\ldots A'_n$, that the open
marked pyramid $OA'_1\ldots A'_n$ contains the marked pyramid
$OA_1\ldots A_n$ (except the point $O$) and the set $OA'_1\ldots
A'_n\setminus OA_1\ldots A_n$ does not contain integer points.

\begin{lemma}\label{last2}
The set $U_n$ is open and nonempty.
Any point of the set $U_n$ defines such two-dimensional continued fraction,
that this fraction contains
the polygon $A_1\ldots A_n$ as a two-dimensional face.
\qed
\end{lemma}

For any convex $k$-gon $P$ in $\r^{n+1}$ for $k\le n+1$
whose two-dimensional plane does not contain the origin,
there exists a $n$-tuple of hyperplanes that divides the two-dimensional
plane containing  $P$ onto connected components, such that
one of these components coincides with $P$.
Further proof of Lemma~\ref{last2} repeats the proof of Lemma~\ref{last}.

\subsubsection{Realizability of faces}

\begin{proposition}
For any $n\ge 2$, any two-dimensional face of List~``$\alpha_n$''
is realizable as a face of some $n$-dimensional continued fraction.
\end{proposition}

\begin{proof}
Since all faces of List~``$\alpha_n$'' ($n \ge 2$)
are triangular or quadrangular (and the corresponding marked pyramids with vertices
at the origin and bases in these faces are completely empty),
Lemmas~\ref{last} and~\ref{last2} can be applied.
\end{proof}

\subsubsection{Nonequivalence of faces}

\begin{lemma}\label{xxxx}
For any  $n\ge 2$, any two different faces of List~``$\alpha_n$''
are integer-linear nonequivalent to each other.
\end{lemma}

\begin{proof}
Lemma~\ref{xxxx} follows directly from Theorem~A (see
List~``M-W'').
\end{proof}

\subsubsection{On polygonal faces of two-dimensional continued fractions}

\begin{corollary}
Let some compact two-dimensional face of some two-dimensional
continued fraction is an integer polygon containing more that three
angles less than the right angle.
Then the integer distance between the origin and the plane containing this face equals one.
\end{corollary}

\begin{proof}
The corollary follows from Lemma~\ref{zapret}.
\end{proof}

\section{Unsolved questions and problems}

In conclusion of this work we outline some arising natural problems here.
First of all let us make the following remark.
By ``classification problems'' for some subset in some other set
through this section we mean the study of the following questions:\\
{\bf a)} which elements of the set are in the subset;\\
{\bf b)} which elements of the set are not in the subset;\\
{\bf c)} which infinite series of elements of the set
are in the subset, how many such series exist;\\
{\bf d)} which infinite series of elements of the set
are not in the subset, how many such series exist;\\
{\bf e)} describe properties of the elements of the subset;\\
{\bf f)} describe properties of the elements of the complement of the
subset in the set;\\
{\bf g)} is the problem of verification weather the given element of the set
is in the subset or not in it algorithmicaly solvable (find the corresponding algorithms);\\
{\bf h)} give the complete list of elements and series of the subset explicitly.

For instance, in this paper we solve the ``classification
problem'' ${\bf h)}$ for the subset of integer-linear or
integer-affine equivalence classes of compact two-dimensional
sails faces of multidimensional continued fractions at the integer
distances from the origin great than one (in the set of
integer-linear/affine equivalence classes of all polygons). Here
the answer to Question~${\bf h)}$ also implies the answers to
Questions~${\bf a)}$, ${\bf b)}$, ${\bf c)}$. Question~${\bf d)}$
becomes meaningless. Also we partially get answers to
Questions~${\bf e)}$ and~${\bf f)}$. Question~${\bf g)}$ is also
closely related to Question~${\bf h)}$ and also was studied by
the author, but it does not appear in the present paper by volume
reasons. The result of this question can be also applied in
algorithms of constructing two-dimensional continuous fractions.

\begin{problem}
Solve the ``classification problems'' for the subset of
integer-linear or integer-affine equivalence classes of compact
three-dimensional $($multi\-di\-men\-sional$)$ sail faces
contained in three-dimensional $($multidimensional$)$ planes at
the integer distance from the origin greater than one $($in the
set of integer-linear/affine equivalence classes of all
polyhedra$)$.
\end{problem}

In connection with the last problem the following question about
marked (compact by definition) pyramids is natural.

\begin{problem}\label{3Dproblem}
Solve all the ``classification problems'' for the subset of
integer-linear or integer-affine equ\-i\-va\-lence clas\-ses of
four-dimensional $($multi\-di\-men\-si\-onal$)$ multistory
com\-pletely empty convex marked pyramids $($in the set of
integer-linear/affine equ\-iv\-a\-lence classes of all convex
marked pyramids of the same dimension$)$.
\end{problem}

The geometrical contents of the next problem is extremely different
from the above two ones.

\begin{problem}
Solve the ``classification problems'' for the subset of
integer-linear or integer-affine equi\-valence clas\-ses of
two/three-dimensional $($multi\-di\-men\-si\-o\-nal$)$ sail faces
contained in the two/three-dimensional
$($multi\-di\-men\-si\-onal$)$ planes at the unit distance from
the origin.
\end{problem}

As a matter of fact this problem can be reduced to the
``classification problems'' for the integer-affine classes of
convex hulls of all integer points in some polygons
$($polyhedra$)$ with bounded number of faces of maximal
dimension. The following question about the polygons and
polyhedra is in its place here.

\begin{problem}
Solve the ``classification problems'' e$)$, f$)$, and h$)$ for the
subset of in\-te\-ger-af\-fine equ\-i\-va\-len\-ce classes of
integer polygons/polyhedra $($in the set of classes of all
polygons/polyhedra$)$.
\end{problem}

The following result on this subject is known.
Denote by $H(\mu)$ the logarithm of the number of integer-affine
equivalence classes of integer polygons of volume
$\mu /2$ in the plane, for some integer $\mu$.

\begin{theorem}{\bf (V.~Arnold~\cite{Arn5}.)}
For sufficiently large $\mu$ the following holds
$$
c_1 \mu ^ {1/3}\le H(\mu) \le c_2 \mu^{1/3}\ln \mu.
$$
\end{theorem}

Further investigations lead to the problems of classifications
of some face arrangements.
We give the simplest intensional example of such problems.

\begin{problem}{\bf(V.~Arnold.)}
Solve the ``classification problems'' for 1-stars of vertex
$($i.e. the union of a vertex and all adjacent edges to this
vertex$)$ for sails of two$($multi$)$-dimensional continued
fractions up to the integer-linear/affine equivalence.
\end{problem}

Here is another problem of this series.

\begin{problem}
Solve the ``classification problems'' for two-tuples of
two-dimensional adjacent faces for the sails of
two$($multi$)$-dimensional continued fractions up to the
integer-linear/affine equivalence.
\end{problem}

The last two problems can be naturally generalized to the case of
more complicated arrangements of faces.

Now we formulate the following problem on so-called {\it stable
integer-affine types} of polyhedra.

\begin{definition}
The integer-affine type of some polyhedron $($polygon$)$ is called
{\it stable in dimension $k$} if for any positive integer $r$
there exists such $k$-dimensional continued fraction that one of
the sails of this fraction contains the face with the given
integer-affine type in the plane at the integer distance equal $r$
to the origin.
\\
The integer-affine type of some polyhedron $($polygon$)$ is called
{\it almost stable in dimension $k$} if for any positive integer
$N$ there exist such $r>N$ and such $k$-dimensional continued
fraction that one of the sails of this fraction contains the face
with the given integer-affine type in the plane at the integer
distance equal $r$ to the origin.
\end{definition}

\begin{problem}
Which integer-affine types of polyhedra are $($almost$)$ stable in
dimension $3$ $($in dimension $k>3)$?
\end{problem}

We can answer on the similar question for the case of polygons.

\begin{corollary}
For any positive integer $k\ge 2$, and any positive integer $a\ge 1$
the integer-affine type of the triangle $(0,0)$, $(a,0)$, $(0,1)$
is stable in dimension $k$ $($see the case $a=6$ on Figure~\ref{qqq}$)$.
\begin{figure}[h]
$$\epsfbox{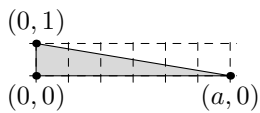}$$
\caption{Stable polygons ($a\ge 1$).}\label{qqq}
\end{figure}
There are no other integer-affine types of integer polygons stable or almost stable
in dimension~$k$.
\qed
\end{corollary}

Even the answer to the following question is unknown to the author.

\begin{problem}
Is it true that the set of all stable in dimension $3$ $($in
dimensiona $k$$)$ integer-affine types of polyhedra coincide with
the set of all almost stable in dimension $3$ $($in dimensiona
$k$$)$ integer-affine types of polyhedra?
\end{problem}

Except the series of problems listed before the problems similar to the
following one are very important and interesting.

\begin{problem}
Do there exist three-dimensional polyhedra that appear
as a faces of sails of $k$-dimensional continued fractions
contained in three-dimensional planes at integer distances
greater than one only for $k>3$, and do not
appear for $k=3$?
In the case of positive answer solve the corresponding ``classification problems''
for them.
\end{problem}

Now we formulate some problems on statistical properties
of sail faces for multidimensional continued fractions.

Denote the set of all integer $(n+1)$-dimensional operators with
real rational distinct eigenvalues by $\Lambda_{n+1}$.
Denote by $\Lambda_{n+1,r}^I$
a subset of $\Lambda_{n+1}$ of operators with the norm not greater than $r$.
Denote by $\Lambda_{n+1,r}^{II}$ a subset of $\Lambda_{n+1}$ of operators with
the norm not greater than $r$ and the square root of the sum of squares of all
characteristic polynomial coefficients not greater than $r$.
(The operator norm here is the square root of the sum of squares of all
its matrix coefficients in some fixed basis.)

Since all eigenvalues of some operator $A$ in the set $\Lambda_{n+1}$
are real and distinct, the number of eigen hyperspaces
for $A$ (in $\r^{n+1}$) equals $n+1$.
The continued fraction is uniquely defined by these hyperspaces.
Since all eigenvalues of $A$ are rational, the sails consist of
finite number of compact faces.
Inversely, if all sails of multidimensional continued fraction consist of
finite number of compact faces, then the continued fraction corresponds to
some operator of the set $\Lambda_{n+1}$.

Let $\Gamma$ be some set of integer-linear types of faces
of $n$-dimensional continued fractions.
By $\sharp_{n+1,r}^I(\Gamma)$ (by $\sharp_{n+1,r}^{II}(\Gamma)$)
we denote the total number of faces with integer-linear type of the set $\Gamma$
for the continued fractions of the set $\Lambda_{n+1,r}^{I}$ ($\Lambda_{n+1,r}^{II}$ respectively).

\begin{problem}{\bf (V.~Arnold.)}
Does there exist a statistics of triangular faces for general
sails of finite multidimensional continued fractions?
Find this statistics in the case of positive answer.
\end{problem}

In other worlds, we have to study the existence of the following limit:
$$
\lim_{r\to \infty}\left(
\frac{\sharp_{n+1,r}^I(\triangle)}{\sharp_{n+1,r}^I(*)}\right)
\quad \left(
 \mbox{or } \lim_{r\to \infty}\left( \frac{\sharp_{n+1,r}^{II}(\triangle)}{\sharp_{n+1,r}^{II}(*)}\right)
\right),
$$
if the limit exists, it is extremely important to find the limit (or even its approximation).
By $\triangle$ we denote the set of all integer-linear types of triangles,
the set of all integer-linear types of faces is denoted by $*$.
Is it true that the limits (for $I$ and $II$) equal to each other?

The similar questions are interesting for the cases of polygons
with $n>3$ vertices, and also for single cases of integer-affine types.
Besides that, the similar questions exist and are interesting for three-dimensional
and multidimensional  polyhedra.

Note that nonexistence of the statistics for some sets of types
does not imply nonexistence of ``relative'' statistics for these sets.

\begin{problem}{\bf (V.~Arnold.)}
Does there exist a ``relative''  statistics of triangular faces
and quadrangular faces for general sails of finite multidimensional continued fractions?
Find this ``relative'' statistics in the case of positive answer.
\end{problem}

As in the previous case we have to study the existence of the following limit (and find it):
$$
\lim_{r\to \infty}\left(
\frac{\sharp_{n+1,r}^I(\triangle)}{\sharp_{n+1,r}^I(\diamondsuit)}\right)
\quad \left(
 \mbox{or } \lim_{r\to \infty}\left( \frac{\sharp_{n+1,r}^{II}(\triangle)}{\sharp_{n+1,r}^{II}(\diamondsuit)}\right)
\right).
$$

Here by $\diamondsuit$ we denote the set of all integer-linear types of quadrangles.
Is it true that the limits (for $I$ and $II$) equal to each other?

\begin{remark}
It is also possible to consider some other exhaustions of $\Lambda_{n+1}$ (except $I$ and $II$)
for calculating the corresponding statistics.
(For more information see, for instance, the work of V.~Arnold~\cite{Arn2}.)
\end{remark}

In the papers~\cite{Arn2},~\cite{Arn4} and the book~\cite{Arn3}
(see problem~1993-11) by V.~Arnold
he discusses notions of statistics for types of sail faces of
multidimensional continued fractions
more detailed and formulates many interesting and actual statistical
problems and conjectures.

For one-dimensional continued fractions some of the conjectures of
V.~Arnold were completely studied by M.~Avdeeva and
V.~Bikovskii~\cite{Avd1} and~\cite{Avd2}. Denote by $\mbox{``}k
\mbox{''}$ a unique integer-linear type of the segment of integer
length $k>0$. In the works~\cite{Avd1} and~\cite{Avd2} for any
$k>0$ M.~Avdeeva and V.~Bikovskii proved the existence and found
the following limits:
$$
\lim_{r\to \infty}\left( \frac{\sharp_{3,r}^I( \mbox{``}k \mbox{''} )}{\sharp_{3,r}^I(*)}\right)
=
\frac{1}{\ln (2)}\ln\left( 1+ \frac{1}{k(k+2)} \right).
$$
and also the authors gave the estimate for the convergence rate of these limits.
It turns out that the limits coincide with the statistics of theorem of Gauss-Kuzmin-L\'evi
(for more information see the works of A.~Wiman~\cite{Wim} and R.~O.~Kuzmin~\cite{Kuz}).

M.~Kontsevich and Yu.~Suhov~\cite{Kon} proved the existence
of an average number of a polyhedron with the prescribed number of integer
points for almost all sails of multidimensional continued fractions
(except some zero Lebesque measure set).
These statistics are not calculated yet, and the methods of their calculation
are not yet developed.
The first statistical data for the periodic two-dimensional continued fractions
is given by the author in the work~\cite{Kar3}.

In conclusion of this section it remains to note that all the problems listed above
can be posed also for the case of sails of periodic algebraic
continued fractions.
We give the following problem as an example.

\begin{problem}{\bf (V.~Arnold.)}
Solve the ``classification problems'' for inte\-ger-affine
equivalence classes of compact two/three$($multi$)$-dimensional
faces of sails of periodic algebraic continued fractions.
\end{problem}

In the last problem it is also useful to study cases of faces contained
in the planes at distances equal/greater than one to the origin.

\begin{remark}
All the statistical questions (similar to the questions
for finite multidimensional continued fractions described above)
can be posed also for the case of periodic multidimensional continued fractions.
\end{remark}

\vspace{10pt} \noindent {\bf Acknowledgements:} \ The author is
grateful to professor V.~I.~Arnold, professor A.~B.~Sossinsky, and
E.~I.~Korkina for constant attention to this work and useful
remarks and discussions, and Universit\'e Paris-Dauphine
--- CEREMADE  for the hospitality and excellent working
conditions.

\end{document}